\documentclass[11pt]{amsart}
\usepackage{amssymb,amsmath}
\usepackage{graphicx}
\usepackage{subcaption}
\usepackage{bm,bbm}
\usepackage{color}
\usepackage{tikz}
\usepackage{pgfplots}
\usepackage{array,blkarray}
\usepackage[hidelinks]{hyperref}
\usepackage{mathtools}
\newtheorem{theorem}{Theorem}

\theoremstyle{definition}

\theoremstyle{lemma}

\theoremstyle{remark}
\newtheorem{remark}[theorem]{Remark}

\numberwithin{theorem}{section}
\numberwithin{equation}{section}
\numberwithin{table}{section}
\numberwithin{figure}{section}
\newcommand{\V}{\ensuremath{\mathcal{V}}}
\newcommand{\Q}{\ensuremath{\mathcal{Q}} }

\def\R{\mathbb{R}}

\definecolor{myBlue1}{RGB}{101,149,239}  
\definecolor{myBlue2}{RGB}{113,104,238} 
\definecolor{myBlue3}{RGB}{30,144,255} 
\definecolor{myGreen1}{RGB}{154,204,50} 
\definecolor{myGreen2}{RGB}{69,169,0} 
\definecolor{myGreen3}{RGB}{154,205,50} 
\definecolor{myGreen4}{RGB}{105,139,34} 
\definecolor{myRed1}{RGB}{210,105,30} 
\definecolor{myRed2}{RGB}{165,42,42} 
\definecolor{myRed3}{RGB}{139,26,26} 
\definecolor{lightgray}{RGB}{175,175,175} 
\definecolor{myLGray}{RGB}{225,225,225} 
\definecolor{mycolor0}{rgb}{0.66,0.66,0.66}
\definecolor{mycolor4}{rgb}{0.00000,0.44700,0.74100}
\definecolor{mycolor1}{rgb}{0.85000,0.32500,0.09800}
\definecolor{mycolor2}{rgb}{0.92900,0.69400,0.12500}
\definecolor{mycolor3}{rgb}{0.67000,0.74700,0.14100}
\definecolor{mycolor5}{rgb}{0.49400,0.18400,0.55600}
\definecolor{mycolor6}{rgb}{0.85000,0.32500,0.09800}%
\DeclareMathOperator{\id}{id}
\DeclareMathOperator{\tr}{tr}

\DeclareMathOperator{\trace}{tr}

\newcommand{\calB}{\ensuremath{\mathcal{B}} }

\newcommand{\calE}{\ensuremath{\mathcal{E}} }

\newcommand{\calJ}{\ensuremath{\mathcal{J}} }
\newcommand{\calK}{\ensuremath{\mathcal{K}} }

\newcommand{\calQ}{\ensuremath{\mathcal{Q}} }
\newcommand{\calR}{\ensuremath{\mathcal{R}} }

\newcommand{\calT}{\ensuremath{\mathcal{T}} }

\def\dx{\,\text{d}x}

\begin{document}
\title[Splitting Schemes for Dynamic Boundary Conditions]{Splitting Schemes for the Semi-linear Wave Equation with Dynamic Boundary Conditions}
\author[]{R.~Altmann$^\dagger$}
\address{${}^{\dagger}$ Institute of Analysis and Numerics, Otto von Guericke University Magdeburg, Universit\"atsplatz 2, 39106 Magdeburg, Germany}
\email{robert.altmann@ovgu.de}
\thanks{Research funded by the support of the Deutsche Forschungsgemeinschaft (DFG, German Research Foundation) through the project 446856041. Moreover, major parts of this work were carried out while the author was affiliated with the Institute of Mathematics and the Centre for Advanced Analytics and Predictive Sciences (CAAPS) at the University of Augsburg. \\
\indent
This paper will appear in Computers and Mathematics with Applications. }
\date{\today}
%
%
\begin{abstract}
This paper introduces novel bulk--surface splitting schemes of first and second order for the wave equation with kinetic and acoustic boundary conditions of semi-linear type. For kinetic boundary conditions, we propose a reinterpretation of the system equations as a coupled system. This means that the bulk and surface dynamics are modeled separately and connected through a coupling constraint. This allows the implementation of splitting schemes, which show first-order convergence in numerical experiments. On the other hand, acoustic boundary conditions naturally separate bulk and surface dynamics. Here, Lie and Strang splitting schemes reach first- and second-order convergence, respectively, as we reveal numerically. 
\end{abstract}
%
\maketitle
%
%
%
\section{Introduction}
We consider wave equations with so-called {\em dynamic boundary conditions} in the semi-linear setting. In contrast to Dirichlet, Neumann, or Robin boundary conditions, such non-trivial boundary conditions do not neglect the momentum of the wave on the boundary. 
The inclusion of boundary dynamics in the model is of particular interest in the field of fluid--structure interaction or acoustic--elastic couplings~\cite{Hip17}. Further applications can be found in connection with Klein--Gordon type equations~\cite{CouFL04}, the stabilization of wave equations via a feed-back law on the boundary~\cite{KomZ90}, or in separation processes in mixtures of two materials~\cite{GarK20}. Moreover, this approach enables to model one part of a coupled problem as a boundary layer, i.e., one wave system may be replaced by a wave-type equation on the boundary~\cite{Lie13,Hip17}. 

In this paper, we consider two different types of boundary conditions: 
{\em Kinetic boundary conditions} are derived in terms of energy balances and constitutive laws and lead to a wave equation in the bulk coupled to a wave equation on the surface. A physical interpretation for the one-dimensional case is given in~\cite{Gol06}. Results on the well-posedness of such models can be found, e.g., in~\cite{Vit13,GraL14}. 
{\em Acoustic boundary conditions} include a second variable for the (small) displacement of the boundary into the domain and were first introduced in~\cite{BeaR74}. Therein, locally reacting boundary conditions were considered, which allow the interpretation of (independent) springs acting in response to the excess pressure in the gas, see also~\cite{GalGG03,Mug06}. 
In general, such boundary conditions model the propagation of sound waves in a fluid at
rest with the second variable modeling the oscillations of the surrounding wall. This can be used, e.g., to model vibrations of the membrane of a bass drum~\cite{Vit15}. For more details and an illustration of the different types of boundary conditions, we refer to~\cite{Hip17}.

Although dynamic boundary conditions are well understood from a theoretical point of view, the literature on numerical methods is rather short. One of the first papers in this direction, considering stationary bulk--surface partial differential equations is~\cite{EllR13}. 
For corresponding parabolic problems, a detailed numerical analysis was introduced in~\cite{KovL17}, see also~\cite{VraS13a}. In both examples, an implicit Euler method is used for the temporal discretization. First approaches in the direction of a bulk--surface splitting were recently introduced in~\cite{AltKZ21,AltZ22ppt_b}. 
For separation processes modeled by the Cahn--Hilliard equation, stable numerical schemes were analyzed in~\cite{Met21,KnoLLM21,BaoZ21}. 
For wave systems with dynamic boundary conditions, the spatial discretization was rigorously analyzed in~\cite{HipHS18,HocL20,HipK20}. 
Moreover, for the temporal discretization, an implicit--explicit variant of the Crank--Nicolson scheme was introduced and analyzed in~\cite{HocL21}. 

This paper is devoted to the construction of novel bulk--surface splitting schemes for wave-type equations with non-standard boundary conditions. Such schemes are of particular interest if the dynamics in the bulk and on the surface have different characteristic length or time scales, cf.~\cite{AltV21}. 
The first step is the derivation of an abstract setting for kinetic as well as acoustic boundary conditions. For kinetic boundary conditions, we consider the equations in the bulk and on the surface as a coupled system of partial differential equations. Together with the necessary coupling condition, this leads to a {\em partial differential-algebraic equation} (PDAE), see~\cite{LamMT13} for an introduction. This formulation then leads to splitting schemes, which decouple the bulk and surface dynamics similarly as in~\cite{AltKZ21} for the parabolic case. To obtain a fully-discrete scheme, this is then combined with the already mentioned bulk--surface finite element method and appropriate time stepping schemes. 
For acoustic boundary conditions, where the bulk and surface dynamics are naturally separated, we proceed contrary and resolve the constraint, i.e., we consider an abstract formulation without an explicit constraint. Based on this formulation, we introduce bulk--surface splitting schemes. The respective combination with an Euler or Crank--Nicolson time discretization then yields splitting schemes of order one and two as we show in numerical experiments.

The paper is organized as follows: in Section~\ref{sect:formulation} we introduce the semi-linear model problems with dynamic boundary conditions and derive an abstract formulation for both types of boundary conditions. 
Splitting schemes for kinetic boundary conditions are then introduced in Section~\ref{sect:schemesKinetic}. For this, we first discuss the spatial discretization and possible time stepping schemes.  
Acoustic boundary conditions are subject of Section~\ref{sect:schemesAcoustic}. Again, we introduce splitting schemes and consider a numerical experiment showing the performance of the newly introduced methods. Finally, we conclude in Section~\ref{sect:conclusion}.
%
%
\section{The Wave Equation with Dynamic Boundary Conditions}
\label{sect:formulation}
As preparation for the construction of splitting methods, this section is devoted to the abstract formulation of the (semi-linear) wave equation with non-trivial boundary conditions. In particular, we introduce possible interpretations of the equations as coupled systems of bulk and surface dynamics. 
%
%
\subsection{Kinetic boundary conditions}
\label{sect:formulation:kinetic}
Let $\Omega$ denote a bounded Lipschitz domain with boundary~$\Gamma\coloneqq\partial\Omega$. As a model problem for the semi-linear wave equation with kinetic boundary conditions, we consider the system 
\begin{subequations}
\label{eq:kineticBC}
\begin{align}
  \ddot u - \Delta u 
  &= f_\Omega(t,u)\qquad \text{in }\Omega, \label{eq:kineticBC:a}\\
  \ddot u - \beta\Delta_\Gamma u + \kappa\, u + \partial_n u 
  &= f_\Gamma(t,u)\qquad \text{on }\Gamma  \label{eq:kineticBC:b}
\end{align}
\end{subequations}
with constants~$\beta, \kappa \ge 0$ and initial conditions~$u(0)=u^0$ and~$\dot u(0)=\dot u^0$. Here, the parameter~$\beta$ in front of the Laplace--Beltrami operator~$\Delta_\Gamma$ (see~\cite[Ch.~16.1]{GilT01}) characterizes whether the boundary conditions are locally reacting ($\beta=0$) or non-local. Throughout this paper, we focus on the latter case, i.e., we assume~$\beta>0$. Finally, $f_\Omega$ and $f_\Gamma$ denote the (sufficiently smooth) nonlinearities in the bulk and on the surface, respectively. Possible extensions of this model problem include advection terms in the bulk and on the surface~\cite{HocL20} 
%
%
or strong damping terms~\cite{HipK20}. 

Problem~\eqref{eq:kineticBC} can be written as an abstract wave equation based on the space $V=\{v\in H^1(\Omega)\ |\ v|_{\Gamma} \in H^1(\Gamma) \}$, cf.~\cite{HipK20,HocL20}. 
This then leads to a variational problem of the form 
\[
  m(\ddot u, v) + d(\dot u, v) + a(u,v) 
  = \langle f(t, u), v\rangle
\]
for test functions~$v\in V$ and corresponding bilinear forms $m$, $d$, $a$, which include bulk as well as surface integrals. In the special case of system~\eqref{eq:kineticBC}, these forms read 
\[
  m(u,v) = \int_\Omega uv \dx + \int_\Gamma uv \dx, \qquad
  a(u,v) = \int_\Omega \nabla u \cdot \nabla v \dx + \int_\Gamma \beta\, \nabla_\Gamma u\cdot \nabla_\Gamma v + \kappa\, uv \dx
\]
together with $d(u,v) = 0$ and $\langle f(t, u), v\rangle = \int_\Omega f_\Omega(t,u)v \dx + \int_\Gamma f_\Gamma(t,u)v \dx$. This abstract formulation allows to use standard methods to prove well-posedness and enables a direct application of a finite element discretization.  
Nevertheless, this formulation is not suitable for the construction of splitting schemes. For parabolic problems with comparable boundary conditions, it has been shown in~\cite{AltKZ21} that a direct splitting approach based on this formulation yields discretization schemes which approximate the solution of a perturbed system. This is caused by the fact that such 'naive' splitting approaches do not include information on the derivatives of the respective variables. 

In order to make the spatial discretization more flexible and to design bulk--surface splitting methods, we introduce an alternative formulation of~\eqref{eq:kineticBC} as a coupled system. We follow the approach of~\cite{Alt19} and introduce, as a first step, the auxiliary variable~$p\coloneqq u|_\Gamma$ on the boundary. With this, system~\eqref{eq:kineticBC} can be written as 
\begin{align*}
	\ddot u - \Delta u 
	&= f_\Omega(t,u)\qquad \text{in }\Omega,\\
	\ddot p - \beta\Delta_\Gamma p + \kappa\, p + \partial_n u 	
	&= f_\Gamma(t,p)\qquad \text{on }\Gamma, \\
	u - p
	&= 0\hspace{1.8cm} \text{on }\Gamma.
\end{align*}
This system consists of two dynamic equations and one constraint, which couples bulk and surface dynamics. In the yet to introduce abstract formulation, we will add the constraint by an additional Lagrange multiplier and introduce the spaces 
\[
  \V \coloneqq H^1(\Omega)\times H^1(\Gamma), \qquad
  \Q \coloneqq H^{-1/2}(\Gamma).
\]
Note that, in contrast to the previous approach, we do not include the fact that $u$ has an $H^1$-trace in the ansatz space but rather have a second variable $p$ with values in~$H^1(\Gamma)$. Hence, we only use standard Sobolev spaces and incorporate information on the trace in form of an additional equation. 

To obtain an abstract formulation, we consider in a second step test functions $(v,q)\in \V$. Multiplying~\eqref{eq:kineticBC:a} by $v$ and~\eqref{eq:kineticBC:b} by $q$ and integrating by parts, we obtain 
\begin{align*}
  (\ddot u, v)_\Omega + (\nabla u, \nabla v)_\Omega - (\partial_n u, v)_\Gamma 
  &= (f_\Omega(t,u), v)_\Omega,\\
  (\ddot p, q)_\Gamma + (\beta\,\nabla_\Gamma p, \nabla_\Gamma q)_\Gamma + (\kappa\, p + \partial_n u, q)_\Gamma 
  &= (f_\Gamma(t,p), q)_\Gamma.
\end{align*}
Here, $(\bullet, \bullet)_\Omega$ and $(\bullet, \bullet)_\Gamma$ denote the standard inner products in~$L^2$ (or the corresponding duality parings) on $\Omega$ and $\Gamma$, respectively. 
For an operator formulation, we introduce~$\calK_\Omega\colon H^1(\Omega) \to [H^1(\Omega)]^*$ and~$\calK_\Gamma\colon H^1(\Gamma) \to [H^1(\Gamma)]^*$ by
\begin{align}
\label{eq:operatorK}
  \langle \calK_\Omega u, v\rangle
  = \int_\Omega \nabla u \cdot \nabla v\dx, \qquad
  \langle \calK_\Gamma p, q\rangle
  = \beta \int_\Gamma \nabla_\Gamma p \cdot \nabla_\Gamma q \dx
  + \kappa \int_\Gamma p \, q \dx.
\end{align}
Moreover, the coupling operator $\calB\colon \V\to\Q^*=H^{1/2}(\Gamma)$ is defined by~$\calB(u,p) \coloneqq p - \trace u$, where~$\trace$ denotes the usual trace operator on~$\Gamma$. With this and a Lagrange multiplier~$\lambda\colon [0,T] \to \Q$, which we insert in place of the normal derivative~$\partial_n u$, we yield the PDAE system 
\begin{subequations}
\label{eq:PDAE:kinetic}
\begin{align}
	\begin{bmatrix} \ddot u \\ \ddot p  \end{bmatrix}
	+ \begin{bmatrix} \calK_\Omega &  \\  & \calK_\Gamma \end{bmatrix}
	\begin{bmatrix} u \\ p  \end{bmatrix}
	+ \calB^*\lambda 
	&= \begin{bmatrix} f_\Omega(t,u) \\ f_\Gamma(t,p) \end{bmatrix} \qquad \text{in } \V^*, \\
	\calB\, \begin{bmatrix} u \\ p  \end{bmatrix} \phantom{i + \calB \lambda} &= \phantom{[]} 0\hspace{5.4em} \text{in } \Q^*.
\end{align}
\end{subequations}
We would like to emphasize that~\eqref{eq:PDAE:kinetic} is indeed equivalent to~\eqref{eq:kineticBC} and that one can show $\lambda = \partial_n u$ under sufficient regularity assumptions on the solution. 
\begin{remark}
The inclusion of additional damping or advection terms does not change the structure of the PDAE in terms of the coupling.  
\end{remark}
\begin{remark}
\label{rem:kinetic:energy}
System~\eqref{eq:PDAE:kinetic} may also be formulated as Hamiltonian (descriptor) system, see~\cite{BeaMXZ18,MehM19} for an introduction. For this, we assume~$f_\Omega=0$, $f_\Gamma=0$ and perform a regularization step in time (known as index reduction in the finite-dimensional setting), i.e., we replace the constraint by its derivative~$\calB(\dot u, \dot p)=\dot p - \trace \dot u =0$. Taking the initial data into account, this yields a formally equivalent system. 

To obtain a first-order formulation, we introduce new variables $w \coloneqq \dot u$, $r \coloneqq \dot p$, which we include in the form $\calK_\Omega \dot u = \calK_\Omega w$ and $\calK_\Gamma \dot p = \calK_\Gamma r$. 
Together with the adapted constraint, this leads to the system 
\begin{align*}
\left[\begin{array}{@{}cc|cc|c}
\calK_\Omega & & & & \\ & \calK_\Gamma & & & \\ \hline 
& & \id & & \\ & & &\id & \\ \hline 
& & & & 0 
\end{array}\right]
\begin{bmatrix} \dot u \\ \dot p \\ \dot w \\ \dot r \\ \dot \lambda \end{bmatrix}
= \left[\begin{array}{@{}cc|cc|c}
& & \calK_\Omega & & \\ 
& &  & \calK_\Gamma & \\ \hline
- \calK_\Omega & & & & \phantom{-}\trace^* \\ 
& - \calK_\Gamma & & & -\id \\ \hline
& & -\trace & \id & 
\end{array}\right]
\begin{bmatrix} u \\ p \\ w \\ r \\ \lambda \end{bmatrix}.  
\end{align*}  
Since~$\calK_\Omega$ and $\calK_\Gamma$ are self-adjoint, the operator matrix on the right is obviously skew-adjoint. Hence, we have a Hamiltonian descriptor system in the sense of~\cite{BeaMXZ18} (see also~\cite{AltMU21b} for the operator case).  
The corresponding energy reads 
\begin{align*}
E(t) 
= \frac 12\ \Big( \|\nabla u\|^2_{L^2(\Omega)} + \beta\, \|\nabla_\Gamma p\|^2_{L^2(\Gamma)} + \kappa\, \|p\|^2_{L^2(\Gamma)} 
+ \|\dot u\|^2_{L^2(\Omega)} + \|\dot p\|^2_{L^2(\Gamma)} \Big) 
\end{align*}
and is preserved over time due to the Hamiltonian structure. 
Note that the first three terms correspond to the potential energy, whereas the last two terms give the kinetic energy of the system. Both parts include contributions from the bulk as well as from the boundary.  		
\end{remark}
We now turn to acoustic boundary conditions, where bulk and surface dynamics are naturally decoupled. 
%
%
\subsection{Acoustic boundary conditions}
\label{sect:formulation:acoustic}
As second example of non-standard boundary conditions, we consider acoustic boundary conditions, which model the propagation of sound waves in a bounded Lipschitz domain~$\Omega$ in combination with oscillations in normal direction on~$\Gamma=\partial\Omega$, cf.~\cite{Hip17}. As in~\cite{HipK20}, we consider the following semi-linear model problem: seek the acoustic velocity potential~$u$ and the (small) displacement of the boundary in normal direction~$\delta$ such that 
\begin{subequations}
\label{eq:acousticBC}
\begin{align} 
  \ddot u - \Delta u  
  &= f_\Omega(t,u)\qquad \text{in }\Omega, \label{eq:acousticBC:a} \\
  \ddot \delta - \beta \Delta_\Gamma \delta + \kappa\, \delta + \dot u
  &= f_\Gamma(t,\delta)\qquad \text{on }\Gamma, \label{eq:acousticBC:b}\\
  \dot \delta 
  &= \partial_n u\hspace{1.35cm} \text{on }\Gamma \label{eq:acousticBC:c}
\end{align}
\end{subequations}
with initial conditions for $u(0)$, $\dot u(0)$ as well as for $\delta(0)$, $\dot \delta(0)$. For the involved constants, we assume~$\beta, \kappa > 0$. Moreover, $f_\Omega$ and $f_\Gamma$ are assumed to be sufficiently smooth as before. 

Similar to the PDAE formulation in the previous subsection, this system contains two dynamic equations and one coupling constraint in~\eqref{eq:acousticBC:c}. In contrast to the kinetic boundary conditions, however, the coupling concerns the normal derivative of~$u$ and is naturally appearing in the system equations. Therefore, we follow the opposite mentality and get rid of the constraint. For this, we first derive the weak formulation. Considering test functions $(v,q)\in \V$ and integrating by parts yields 
\begin{align*}
  (\ddot u, v)_\Omega + (\nabla u, \nabla v)_\Omega - (\partial_n u, v)_\Gamma 
  &= (f_\Omega(t,u), v)_\Omega, \\
  (\ddot \delta, q)_\Gamma + (\beta\,\nabla_\Gamma \delta, \nabla_\Gamma q)_\Gamma + (\kappa\, \delta + \dot u, q)_\Gamma 
  &= (f_\Gamma(t,\delta), q)_\Gamma.
\end{align*}
Replacing the normal derivative of $u$ in the first equation by the constraint~$\dot \delta = \partial_n u$ and using the operators introduced in~\eqref{eq:operatorK}, we obtain the operator equation 
\begin{align}
\label{eq:PDAE:acoustic}
	\begin{bmatrix} \ddot u \\ \ddot \delta  \end{bmatrix}
	+ \begin{bmatrix} & -\tr^* \\ \trace & \end{bmatrix}
	\begin{bmatrix} \dot u \\ \dot \delta \end{bmatrix}
	+ \begin{bmatrix} \calK_\Omega &  \\  & \calK_\Gamma \end{bmatrix}
	\begin{bmatrix} u \\ \delta \end{bmatrix}
	&= \begin{bmatrix} f_\Omega(t,u) \\ f_\Gamma(t,\delta) \end{bmatrix} \qquad \text{in } \V^*.  
\end{align}
Note that, in this formulation, the coupling of $u$ and $\delta$ appears through the terms involving the first derivatives. Further note that this formulation is equivalent to the weak formulation derived in~\cite{HipK20}, where it is written in terms of bilinear forms. 
\begin{remark}
\label{rem:acoustic:energy}
We comment on the Hamiltonian structure of the system, for which we consider the first-order formulation with new variables~$w \coloneqq \dot u$ and $\zeta \coloneqq \dot \delta$. Assuming~$f_\Omega=0$, $f_\Gamma=0$, this yields 
\begin{align*}
\begin{bmatrix} \dot u \\ \dot \delta \\ \dot w \\ \dot \zeta \end{bmatrix}
= \left[\begin{array}{@{}cc|cc@{\ }}
& & \phantom{-}\id &  \\ 
& & & \id \\ \hline
-\id & & \phantom{-}0 & \trace^* \\ 
& -\id & -\trace & 0
\end{array}\right]
\left[\begin{array}{@{}cc|cc@{\,}}
\calK_\Omega & & & \\ & \calK_\Gamma & & \\ \hline & & \id & \\ & & & \id
\end{array}\right]
\begin{bmatrix} u \\ \delta \\ w \\ \zeta \end{bmatrix}, 
\end{align*}
which is a Hamiltonian system of the form~$\calE \dot z = (\calJ-\calR)\calQ z$, see again~\cite{BeaMXZ18}, with $\calE=\id$, $\calR=0$, and $\calJ$ being skew-adjoint. In this form, the energy is defined via $E(t) = \frac 12 \langle \calE^*\calQ z,z\rangle$, which leads to 
\begin{align*}
E(t) 
= \frac 12\ \Big( \|\nabla u\|^2_{L^2(\Omega)} + \beta\, \|\nabla_\Gamma \delta\|^2_{L^2(\Gamma)} + \kappa\, \|\delta\|^2_{L^2(\Gamma)} 
+ \|\dot u\|^2_{L^2(\Omega)} + \|\dot \delta\|^2_{L^2(\Gamma)} \Big). 
\end{align*}
As before, this energy is preserved due to the Hamiltonian structure of the system. 
\end{remark}
%
Based on the abstract formulations presented in this section, we now turn to the construction of splitting methods, which decouple bulk and surface dynamics. We start with kinetic boundary conditions before we consider acoustic boundary conditions in the ensuing section. 
%
%
\section{Splitting Schemes for Kinetic Boundary Conditions}\label{sect:schemesKinetic}
Due to its saddle point structure, a spatial discretization of the PDAE~\eqref{eq:PDAE:kinetic} by finite elements leads to a semi-explicit differential-algebraic equation (DAE) of index 3. Based on a regularization of the resulting semi-discrete system, we introduce a splitting of bulk and surface dynamics. Together with suitable time stepping schemes, this then leads to fully-discrete splitting schemes. Unfortunately, numerical experiments indicate that the resulting schemes are at most of order one (in time). 
%
%
\subsection{Spatial discretization}\label{sect:schemesKinetic:FEM}
For the discretization in space, we consider bulk--surface finite elements, which we only shortly discuss. More details can be found in~\cite{EllR13} as well as in~\cite{KovL17}. 

The spatial domain~$\Omega$ is approximated by a quasi-uniform family of meshes~$\calT_h$ with maximal mesh width $h$. Since the boundary of $\Omega$ may be curved, the union of all elements of~$\calT_h$ defines a polyhedral domain $\Omega_h$ (which may differ from $\Omega$) with boundary $\Gamma_h$. Throughout this paper, we assume that the vertices of~$\Gamma_h$ are part of~$\Gamma$, cf.~\cite{EllR13}. 

Given a mesh~$\calT_h$, we consider the standard $P_1$-finite element space. Note that this yields a nonconforming approximation of $H^1(\Omega)$ if~$\Omega_h\neq \Omega$. A suitable lift operator is introduced in~\cite{Dzi88}. A corresponding basis is given by the usual globally continuous and piecewise linear nodal basis functions. For the discretization of~$p$, we use the mesh~$\calT_h$ restricted to the boundary of~$\Omega_h$. 
This choice leads to the mass matrices $M_\Omega\in \R^{N_\Omega,N_\Omega}$ and $M_\Gamma\in \R^{N_\Gamma,N_\Gamma}$ as discrete versions of the respective $L^2$-inner products and the stiffness matrices $A_\Omega\in \R^{N_\Omega,N_\Omega}$ and $A_\Gamma\in \R^{N_\Gamma,N_\Gamma}$ as discrete versions of the operators~$\calK_\Omega$ and $\calK_\Gamma$, respectively. 
The discrete version of the coupling operator~$\calB$ is denoted by $B\in \R^{N_\Gamma, N_\Omega+N_\Gamma}$ and has full row-rank. Moreover, assuming the a specific ordering of the basis functions such that the last nodes are on the surface, we get $B=[\, 0\ \ M_\Gamma\ -\!M_\Gamma]$. 

The resulting semi-discrete system, where we seek~$u\colon[0,T] \to \R^{N_\Omega}$, $p\colon[0,T] \to \R^{N_\Gamma}$, and the Lagrange multiplier~$\lambda\colon[0,T] \to \R^{N_\Gamma}$, reads 
\begin{subequations}
\label{eq:semidisc:kinetic}
\begin{align}
	\begin{bmatrix} M_\Omega &  \\  & M_\Gamma \end{bmatrix}
	\begin{bmatrix} \ddot u \\ \ddot p  \end{bmatrix}
	+ \begin{bmatrix} A_\Omega &  \\  & A_\Gamma \end{bmatrix}
	\begin{bmatrix} u \\ p  \end{bmatrix}
	+ B^T \lambda 
	&= \begin{bmatrix} f_\Omega(t,u) \\ f_\Gamma(t,p) \end{bmatrix}, \label{eq:semidisc:kinetic:a} \\
	\begin{bmatrix} 0 & M_\Gamma \end{bmatrix} u - M_\Gamma p 
	&= 0. \label{eq:semidisc:kinetic:b}
\end{align}
\end{subequations}
Due to the construction, we assume that $M_\Omega, M_\Gamma$ are symmetric and positive definite, whereas~$A_\Omega, A_\Gamma$ are symmetric and semi-positive definite. With the full-rank property of~$B$, this implies that system~\eqref{eq:semidisc:kinetic} indeed equals a DAE of index 3, cf.~\cite[Ch.~VII.1]{HaiW96}. 
%
%
\subsection{Temporal discretization schemes for wave-type equations}\label{sect:schemesKinetic:time}
As final preparation for the construction of fully-discrete splitting schemes, we recall time stepping schemes for systems of the form $M \ddot u + D \dot u + A u = f$ on a time interval~$[0,T]$. For simplicity, we assume equidistant time steps of size~$\tau$, leading to discrete time points~$t^n\coloneqq n\tau$. We first consider the {\em implicit Euler scheme} applied to the corresponding first-order formulation with ~$w\coloneqq\dot u$. This gives $M w^{n+1} + \tau D w^{n+1} + \tau A u^{n+1} = M w^{n} + \tau f^{n+1}$. Then, replacing~$u^{n+1}$ by the equation $u^{n+1} = u^n + \tau\, w^{n+1}$, we obtain the discretization scheme 
\begin{subequations}
\label{eq:implicitEuler}
\begin{align}
  \big( M + \tau D + \tau^2 A \big)\, w^{n+1}  
  &= M w^n - \tau A u^n + \tau f^{n+1}, \\
  u^{n+1}
  &= u^n + \tau\, w^{n+1}.
\end{align}
\end{subequations}
In the linear case, we have $f^{n+1} = f(t^{n+1})$, whereas in the semi-linear case, one may include the nonlinearity explicitly, i.e., $f^{n+1} = f(t^{n+1},u^n)$. Note that the two equations~\eqref{eq:implicitEuler} can be computed sequentially, i.e., in two steps. 

In a similar way, one may also consider the trapezoidal rule in such a two-step formulation. 
%
%
%
Finally, we mention the {\em Crank--Nicolson scheme}, which can be written in a three-step formulation. A direct application of Crank--Nicolson to the first-order formulation and replacing $u^{n+1}$ as before yields 
\begin{align*}
	\big( M + \tfrac\tau2 D + \tfrac{\tau^2}4 A \big) (w^n + w^{n+1}) 
	&= 2 M w^{n} - \tau A u^n + \tfrac\tau2\, (f^n + f^{n+1}), \\
	u^{n+1} 
	&= u^n + \tfrac\tau2\, (w^n + w^{n+1}). 
\end{align*}
Then, introducing $w^{n+1/2} \coloneqq \frac 12 (w^n + w^{n+1})$
%
%
and changing the incorporation of the right-hand side slightly, namely replacing the trapezoidal rule by a left/right rectangle rule, leads to a Crank--Nicolson scheme with an explicit implementation of the nonlinearity, cf.~\cite{HocL21}. The resulting implicit--explicit scheme reads 
\begin{subequations}
\label{eq:CrankNicolson}
\begin{align}
  \big( M + \tfrac{\tau}{2} D + \tfrac{\tau^2}{4} A \big)\, w^{n+1/2} 
  &= M w^n - \tfrac{\tau}{2} A u^n + \tfrac{\tau}{2} f^n, \\
  u^{n+1}
  &= u^n + \tau\, w^{n+1/2}, \\ 
  M w^{n+1} 
  &= 2M w^{n+1/2} - M w^n + \tfrac\tau2\, ( f^{n+1}-f^n ). 
\end{align}
\end{subequations}
Here, we write $f^{n} = f(t^n,u^n)$ and $f^{n+1} = f(t^{n+1},u^{n+1})$, i.e., $f^{n+1}$ includes the already updated approximation $u^{n+1}$. 
%
%
\subsection{Semi-discrete system and bulk--surface splitting}
\label{sect:schemesKinetic:semiDiscrete}
In order to construct bulk--surface splitting schemes, we need a decomposition of the mass and stiffness matrices. Using the already mentioned convention that the last degrees of freedom correspond to the surface, we obtain
\begin{align}
\label{eq:decomposition}
  u 
  = \begin{bmatrix} u_1 \\ u_2 \end{bmatrix}, \qquad 
  M_\Omega 
  = \begin{bmatrix} M_{11} & M_{12} \\ M_{21} & M_{22} \end{bmatrix}, \qquad 
  A_\Omega 
  = \begin{bmatrix} A_{11} & A_{12} \\ A_{21} & A_{22} \end{bmatrix}
\end{align}
with $u_1(t)\in\R^{N_\Omega-N_\Gamma}$, $u_2(t)\in\R^{N_\Gamma}$ and $M_{ij}$, $A_{ij}$ of corresponding dimension. In the same way, the right-hand side~$f_\Omega$ is decomposed into $f_1$ and~$f_2$.

For the construction of splitting methods, we introduce two subsystems which decouple the bulk and surface dynamics. For the sake of clarity, we consider the first time interval~$[0,\tau]$ only with $\tau$ being the step size. 
As initial data, we assume given values~$u(0) = u^0$ (also defining $u_1^0, u_2^0$), $\dot u(0) = \dot u^0$, $p(0) = p^0 = u_2^0$, and~$\dot p(0) = \dot p^0$. Moreover, we assume that we also have an approximation of the second derivative, which we denote by~$\ddot p^0$. 
Succeeding the approach for parabolic systems~\cite{AltKZ21}, we propose the following decomposition: 

{\em Subsystem 1 (bulk problem)}: 
We consider the dynamics in the bulk with Dirichlet boundary conditions given by the current value of~$p$ and its temporal derivatives, cf.~\cite{Alt15}. 
%
%
With the decomposition of the mass and stiffness matrices introduced in~\eqref{eq:decomposition} and the special structure of~$B$, this then leads to 
\begin{align}
\label{eq:kinetic:bulk}
	M_{11}{\ddot u}_1 + A_{11} u_1 
	= f_{1}(u) - M_{12}\ddot p^0 - A_{12} p^0 
\end{align}
with initial conditions $u_1(0) = u_1^0$ and~$\dot u_1(0) = \dot u_1^0$. Since we do not update~$p$ within this first subsystem, we get~$u_2 \equiv p^0$ due to the constraint. 
Solving~\eqref{eq:kinetic:bulk} on~$[0,\tau]$, we end up with an approximation~$u_1(\tau)$ as well as corresponding derivatives. 

{\em Subsystem 2 (boundary problem)}: 
The second subsystem is a pure boundary problem and reads 
\begin{align}
\label{eq:kinetic:boundary}
	M_\Gamma \ddot p + A_\Gamma p 
	= f_\Gamma(p) + f_{2}(u) 
	- M_{22} \ddot p^0 - A_{22} p^0
	- M_{21} {\ddot u}_1(\tau) - A_{21} u_1(\tau)
\end{align}
with initial conditions~$p(0) = p^0$ and~$\dot p(0) = \dot p^0$. Here, we fix the values coming from the bulk (i.e., $u_1$), as they were already updated. 
From~\eqref{eq:kinetic:boundary} we obtain the new values of $p$ at time~$t=\tau$, which are then inserted in the bulk problem~\eqref{eq:kinetic:bulk} of the subsequent time step. 

The presented decomposition already indicates the possibility of a bulk--surface splitting method. For practical simulations, however, an additional temporal discretization is necessary. 
%
%
\subsection{Fully-discrete splitting schemes}\label{sect:schemesKinetic:fullyDiscrete}
We propose two fully-discrete splitting schemes. First, we consider a Lie splitting, where the two subsystems are solved by an implicit Euler discretization. Second, we introduce a symmetric splitting approach, where we implement the implicit--explicit Crank--Nicolson discretization in each subsystem. 
%
%
\subsubsection{Lie splitting with implicit Euler discretization}\label{sect:schemesKinetic:fullyDiscrete:Euler}
For the Lie splitting approach, where we solve the bulk and boundary problems sequentially, we can expect at most first-order convergence. Hence, we implement a simple Euler discretization for both subsystems. More precisely, we propose to apply the implicit Euler method with an explicit treatment of the possible nonlinearity. 
Hence, we apply~\eqref{eq:implicitEuler} to the two subsystems~\eqref{eq:kinetic:bulk} and~\eqref{eq:kinetic:boundary}. 

One step of Lie splitting then reads as follows: 
Given initial data $u_1^n$, $w_1^n\coloneqq\dot u_1^n$, $p^n$, $r^n\coloneqq \dot p^n$, and~$\ddot p^n$, solve the fully-discrete bulk problem 
\begin{align*}
  \big( M_{11} + \tau^2 A_{11} \big)\, w_1^{n+1}  
  &= M_{11} w_1^n - \tau A_{11} u_1^n + \tau\, (f_1^{n+1} - M_{12}\ddot p^n - A_{12} p^n), \\
  u_1^{n+1}
  &= u_1^n + \tau\, w_1^{n+1}. 
\end{align*}
Recall that, in the semi-linear case, $f_1^{n+1}$ should be understood as $f_1(t^{n+1},u^n)$. Moreover, we obtain an approximation of the second derivative by $\ddot u_1^{n+1} \coloneqq \frac1\tau(w_1^{n+1}-w_1^n)$. 
%
Introducing 
\[
  g^{n+1} 
  \coloneqq f_\Gamma^{n+1} + f_{2}^{n+1} 
  - M_{22} \ddot p^n - A_{22} p^n
  - M_{21} \ddot u_1^{n+1} - A_{21} u_1^{n+1}, 
\] 
the Euler scheme applied to the boundary subsystem reads 
\begin{align*}
  \big( M_\Gamma + \tau^2 A_\Gamma \big)\, r^{n+1}  
  &= M_\Gamma r^n - \tau A_\Gamma p^n + \tau\, g^{n+1}, \\
  p^{n+1}
  &= p^n + \tau\, r^{n+1}.
\end{align*}
As approximation of the second derivative of $p$ we set~$\ddot p^{n+1}\coloneqq\frac1\tau(r^{n+1}-r^n)$. 
Corresponding numerical experiments are subject of Section~\ref{sect:schemesKinetic:numerics}. 
%
%
\subsubsection{Strang splitting with Crank--Nicolson discretization}\label{sect:schemesKinetic:fullyDiscrete:CN}
In this second approach, we consider a symmetric splitting. This means that we first solve the bulk problem~\eqref{eq:kinetic:bulk} on the first half of the respective interval. Then, the boundary problem~\eqref{eq:kinetic:boundary} is solved on the entire interval before we close with the remaining part of the bulk problem. Here, each subsystem is discretized by the implicit--explicit Crank--Nicolson scheme given in~\eqref{eq:CrankNicolson}.

The resulting Strang splitting scheme reads as follows: 
Given initial data $u_1^n$, $w_1^n$, $p^n$, $r^n$, and~$\ddot p^n$ as before, solve 
\begin{align*}
  \big( M_{11} + \tfrac{\tau^2}{16} A_{11} \big)\, w_1^{n+1/4}  
  &= M_{11} w_1^n - \tfrac{\tau}{4} A_{11} u_1^n + \tfrac{\tau}{4}\, (f_1^{n} - M_{12}\ddot p^n - A_{12} p^n), \\
  u_1^{n+1/2}
  &= u_1^n + \tfrac\tau2\, w_1^{n+1/4}, \\ 
  M_{11} w_1^{n+1/2}
  &= 2M_{11} w_1^{n+1/4} - M_{11} w_1^{n} + \tfrac{\tau}{4}\, (f_1^{n+1/2}-f_1^n).
\end{align*}
Recall that, in the semi-linear setting, $f_1^{n+1/2}$ involves the updated state~$u_1^{n+1/2}$ (but still $u_2^n$). For the approximation of the second derivative, we set $\ddot u_1^{n+1/2} \coloneqq \frac2\tau\, (w_1^{n+1/2}-w_1^n)$. 

In the second step, we solve the boundary problem on the entire time interval. 
Introducing 
\[
  h^{n}	
  \coloneqq - M_{22} \ddot p^n - A_{22} p^n 
  - M_{21} \ddot u_1^{n+1/2} - A_{21} u_1^{n+1/2},
\] 
we compute 
\begin{align*}
	\big( M_\Gamma + \tfrac{\tau^2}{4} A_\Gamma \big)\, r^{n+1/2} 
	&= M_\Gamma r^n - \tfrac{\tau}{2} A_\Gamma p^n + \tfrac\tau2\, ( f_\Gamma^{n} + f_{2}^{n} + h^n ), \\
	p^{n+1}
	&= p^n + \tau\, r^{n+1/2}, \\ 
	M_\Gamma r^{n+1} 
	&= 2M_\Gamma r^{n+1/2} - M_\Gamma r^n + \tfrac{\tau}{2}\, ( f_\Gamma^{n+1} + f_2^{n+1} - f_\Gamma^n - f_2^n ). 
\end{align*}
As before, we define~$\ddot p^{n+1} \coloneqq \frac1\tau\, (r^{n+1}-r^n)$.  
%
Finally, we solve the bulk problem on the second half of the time interval. Hence, we compute 
\begin{align*}
  \big( M_{11} + \tfrac{\tau^2}{16} A_{11} \big)\, w_1^{n+3/4}  
  &= M_{11} w_1^{n+1/2} - \tfrac{\tau}{4} A_{11} u_1^{n+1/2} + \tfrac\tau4\, (f_1^{n+1/2} - M_{12}\ddot p^{n+1} - A_{12} p^{n+1}), \\
  u_1^{n+1}
  &= u_1^{n+1/2} + \tfrac\tau2\, w_1^{n+3/4}, \\ 
  M_{11} w_1^{n+1}
  &= 2M_{11} w_1^{n+3/4} - M_{11} w_1^{n+1/2} + \tfrac\tau4\, (f_1^{n+1}-f_1^{n+1/2}).
\end{align*}
We now turn to numerical experiments in order to explore the potential of the proposed splitting schemes. 
%
%
\subsection{Numerical experiments}\label{sect:schemesKinetic:numerics}
This part is devoted to the experimental investigation of the convergence orders of the proposed splitting methods. Moreover, we are interested in possible $h$-dependencies, i.e., we ask ourselves in which way the convergence of the splitting scheme depends on the spatial discretization parameter. 

The spatial domain is given by the unit disc $\Omega = \{x\in \R^2\,|\, x^2_1 + x^2_2 \leq 1\}$ and meshes are generated by {\em DistMesh}~\cite{PerS04}. In all experiments, we consider the error between the numerical solution and a reference solution, which is obtained by a Crank--Nicolson discretization of the full system (without a splitting) with step size~$\tau_\text{ref}=2^{-12}$. Since we focus on the error caused by the splitting and the temporal discretization, the numerical solution acts on the same spatial mesh as the reference solution. 
The errors are measured in the norms corresponding to the spaces 
\[
	L^\infty(L^2(\Omega)) \coloneqq L^\infty(0,T;L^2(\Omega)), \qquad
	L^\infty(H^1(\Omega)) \coloneqq L^\infty(0,T;H^1(\Omega))	
\]
for $u$ as well as 
\[
	L^\infty(L^2(\Gamma)) \coloneqq L^\infty(0,T;L^2(\Gamma)), \qquad
	L^\infty(H^1(\Gamma)) \coloneqq L^\infty(0,T;H^1(\Gamma))	
\]
for $p$.
\subsubsection*{Initial data and parameters}
Following \cite[Sect 8.1]{HipK20}, we set as initial data
\[
	u^0(x,y) = \exp(-20((x-1)^2+y^2)), \qquad
	\dot u^0(x,y) = 0. 
\]
The values for $p^0$ and $\dot p^0$ are chosen in a consistent manner, i.e., $p^0$ equals $u^0$ restricted to the boundary and~$\dot p^0=0$. For the parameters in~\eqref{eq:kineticBC}, we set 
\[
	\beta = 1, \qquad
	\kappa = 1, \qquad 
	T = 1.
\]
In order to analyze the dependence on the mesh size, we perform the upcoming experiments on two different triangulations with mesh sizes 
%
\[
	h_\text{coarse} \approx 0.094, \qquad 
	h_\text{fine} \approx 0.025.	
\]
Moreover, in order to track the convergence w.r.t.~the step size~$\tau$, we consider step sizes $\tau = 2^{-4}, \dots, 2^{-11}$. 
\subsubsection*{Order of convergence in the linear case}
\begin{figure}
	\includegraphics[width=.3\textwidth]{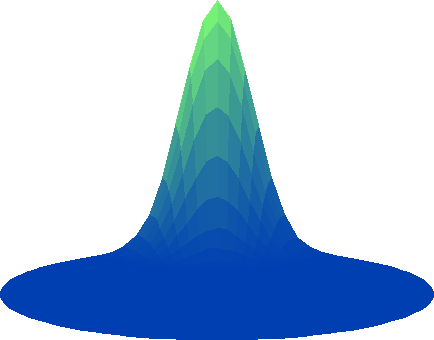}\hspace{0.5em}
	\includegraphics[width=.3\textwidth]{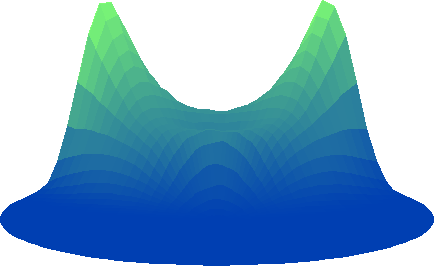}\hspace{0.5em}
	\includegraphics[width=.3\textwidth]{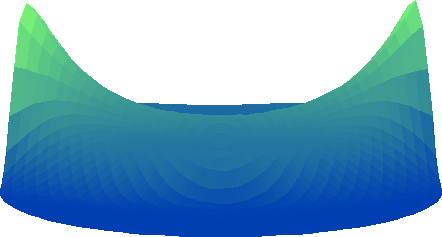}\\[0.5em]
	\includegraphics[width=.3\textwidth,height=.3\textwidth]{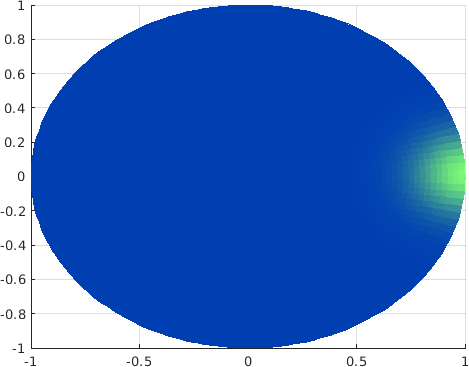}\hspace{0.5em}
	\includegraphics[width=.3\textwidth,height=.3\textwidth]{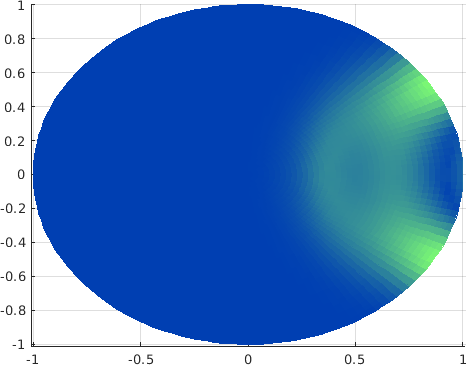}\hspace{0.5em}
	\includegraphics[width=.3\textwidth,height=.3\textwidth]{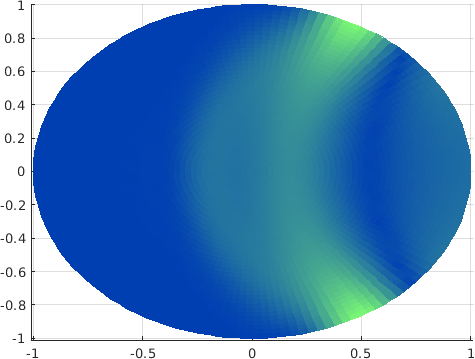}
	\caption{Snapshots of the solution $u$ at times $t=0$ (initial data, left), $t=0.5$ (middle), and $t=1.0$ (right) from two different perspectives. The lower pictures show the top view, whereas the upper plots show the approaching wave from the point $(-1,0)$.}	
	\label{fig:kinetic:solution}	
\end{figure}

\begin{figure} 
%
%
\begin{tikzpicture}
\begin{axis}[%
width=2.3in,
height=1.8in,
scale only axis,
xmode=log,
xmin=0.00045,
xmax=0.06,
xminorticks=true,
xlabel style={font=\color{white!15!black}},
xlabel={step size~$\tau$},
ymode=log,
ymin=2e-04,
ymax=6e-01,
yminorticks=true,
ylabel style={font=\color{white!15!black}},
axis background/.style={fill=white},
title style={font=\bfseries},
legend columns = 2,
legend style={at={(1.89,1.02)}, anchor=south east, legend cell align=left, align=left, draw=white!15!black}
]


\addplot [color=mycolor2, line width=1.5]
table[row sep=crcr]{%
	0.0625	0.077528392\\
	0.03125	0.050674864\\
	0.015625	0.031064639\\
	0.0078125	0.018125693\\
	0.00390625	0.010137677\\
	0.001953125	0.0054699528\\
	0.0009765625	0.0028730025\\
	0.00048828125	0.0014818236\\
};
\addlegendentry{Lie--Euler with $h_\text{fine}$\qquad} 

\addplot [color=mycolor2, line width=1.5, dashed]
table[row sep=crcr]{%
	0.0625	0.070229358\\
	0.03125	0.045449164\\
	0.015625	0.027710399\\
	0.0078125	0.016189398\\
	0.00390625	0.009104677\\
	0.001953125	0.0049442607\\
	0.0009765625	0.0026124602\\
	0.00048828125	0.0013533028\\
};
\addlegendentry{Lie--Euler with $h_\text{coarse}$} 


\addplot [color=mycolor5, line width=1.5]
table[row sep=crcr]{%
	0.0625	0.21308087\\
	0.03125	0.08949316\\
	0.015625	0.04033554\\
	0.0078125	0.019140018\\
	0.00390625	0.0093245097\\
	0.001953125	0.0046023871\\
	0.0009765625	0.0022864012\\
	0.00048828125	0.0011395272\\
};
\addlegendentry{Strang--CN with $h_\text{fine}$\qquad} 

\addplot [color=mycolor5, line width=1.5, dashed]
table[row sep=crcr]{%
	0.0625	0.035064267\\
	0.03125	0.016340245\\
	0.015625	0.0079197404\\
	0.0078125	0.0039076193\\
	0.00390625	0.0019430074\\
	0.001953125	0.00096911109\\
	0.0009765625	0.00048399765\\
	0.00048828125	0.00024186623\\
};
\addlegendentry{Strang--CN with $h_\text{coarse}$} 


\addplot [color=gray, dotted, line width=1.0pt]
table[row sep=crcr]{%
	0.0625	0.6\\
	0.0000625	0.0006\\
};


\end{axis}

%
%
\begin{axis}[%
width=2.3in,
height=1.8in,
at={(2.45in,0.in)},
scale only axis,
xmode=log,
xmin=0.00045,
xmax=0.06,
xminorticks=true,
xlabel style={font=\color{white!15!black}},
xlabel={step size~$\tau$},
ymode=log,
ymin=2e-04,
ymax=6e-01,
yminorticks=true,
yticklabel pos=right,
ylabel style={font=\color{white!15!black}},
axis background/.style={fill=white},
title style={font=\bfseries},
]


\addplot [color=mycolor2, line width=1.5]
table[row sep=crcr]{%
	0.0625	0.16724683\\
	0.03125	0.11755939\\
	0.015625	0.074564637\\
	0.0078125	0.043280888\\
	0.00390625	0.023575486\\
	0.001953125	0.012346029\\
	0.0009765625	0.0063236752\\
	0.00048828125	0.0032010273\\
};

\addplot [color=mycolor2, line width=1.5, dashed]
table[row sep=crcr]{%
	0.0625	0.16417445\\
	0.03125	0.11558231\\
	0.015625	0.073567554\\
	0.0078125	0.042885518\\
	0.00390625	0.023446349\\
	0.001953125	0.012309193\\
	0.0009765625	0.0063140882\\
	0.00048828125	0.0031987237\\
};


\addplot [color=mycolor5, line width=1.5]
table[row sep=crcr]{%
	0.0625	0.59597816\\
	0.03125	0.23770669\\
	0.015625	0.10492196\\
	0.0078125	0.049329082\\
	0.00390625	0.023927432\\
	0.001953125	0.011785139\\
	0.0009765625	0.0058486552\\
	0.00048828125	0.0029134542\\
};

\addplot [color=mycolor5, line width=1.5, dashed]
table[row sep=crcr]{%
	0.0625	0.10224734\\
	0.03125	0.046022227\\
	0.015625	0.02210816\\
	0.0078125	0.010891533\\
	0.00390625	0.0054134433\\
	0.001953125	0.0026996974\\
	0.0009765625	0.0013482417\\
	0.00048828125	0.00067375568\\
};


\addplot [color=gray, dotted, line width=1.0pt]
table[row sep=crcr]{%
	0.0625	0.6\\
	0.0000625	0.0006\\
};


\end{axis}
\end{tikzpicture}%
	\caption{Convergence history for the linear case with kinetic boundary conditions. Plots show the $L^\infty(L^2(\Omega))$ error of $u$ (left) and the $L^\infty(L^2(\Gamma))$ error of~$p$ (right) for the two splitting schemes and two spatial mesh sizes $h_\text{fine}$ (solid) and $h_\text{coarse}$ (dashed). The dotted line indicates order 1.}
	\label{fig:kinetic:linear}
\end{figure}
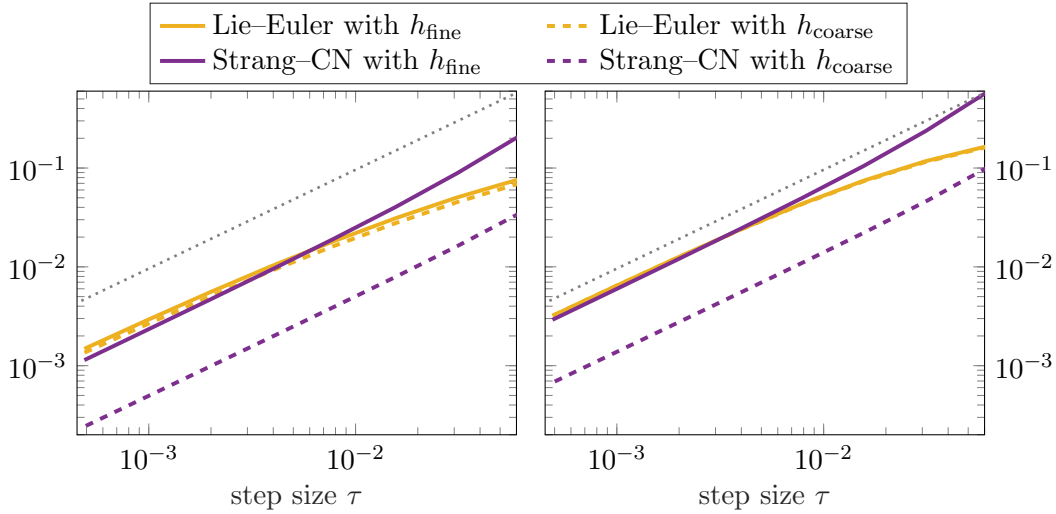

In our first experiment, we consider the {\em linear} case with $f_\Omega=0$, $f_\Gamma=0$. An illustration of the solution evolving over time is shown in Figure~\ref{fig:kinetic:solution}. 
Regarding the convergence w.r.t.~the $L^\infty(L^2(\Omega))$ and $L^\infty(L^2(\Gamma))$-norms for $u$ and $p$, respectively, we refer to Figure~\ref{fig:kinetic:linear}.
%
%
Therein, one can see that the Lie--Euler splitting introduced in~Section~\ref{sect:schemesKinetic:fullyDiscrete:Euler} converges (almost) with order 1 in $u$ and in $p$. More precisely, we obtain the rates $0.92$ for $u$ and $0.96$ for $p$. Note that the dashed and solid yellow lines are on top of each other, meaning that the convergence is independent of the spatial discretization parameter. Besides this stability w.r.t.~the spatial mesh, the splitting scheme is also stable with regard to very coarse time steps. 
\begin{remark} 
If we replace the implicit Euler discretization by the Crank--Nicolson method, then the results are twofold. First, we obtain full first-order convergence, also in the respective $L^\infty(H^1(\Omega))$ and $L^\infty(H^1(\Gamma))$-norms. On the other hand, there occurs a clearly observable $h$-dependence. By this we mean that the error curves move upwards for refinements of the spatial mesh. As a consequence, there cannot be a convergence result for this combination in terms of~$\tau$, which is independent of $h$. Convergence may only be proven in the presence of a CFL-type condition ensuring stability. 
\end{remark}
Unfortunately, the combination of Strang splitting and the Crank--Nicolson discretization as introduced in Section~\ref{sect:schemesKinetic:fullyDiscrete:CN} does not improve the convergence. This is in line with the corresponding parabolic case analyzed in~\cite{AltKZ21}, where the maximal order of two is neither reached. Moreover, the errors are $h$-dependent. This can be seen in Figure~\ref{fig:kinetic:linear}, where the purple solid line (corresponding to the refined spatial discretization) shows an error approximately 5 times larger than on the coarse mesh.  

To summarize, the fully-discrete Lie--Euler splitting shows promising results with first-order convergence independent of the spatial mesh size. The Strang--CN scheme, on the other hand, neither reaches second-order convergence nor is independent of $h$. 
\subsubsection*{Conservation of energy}
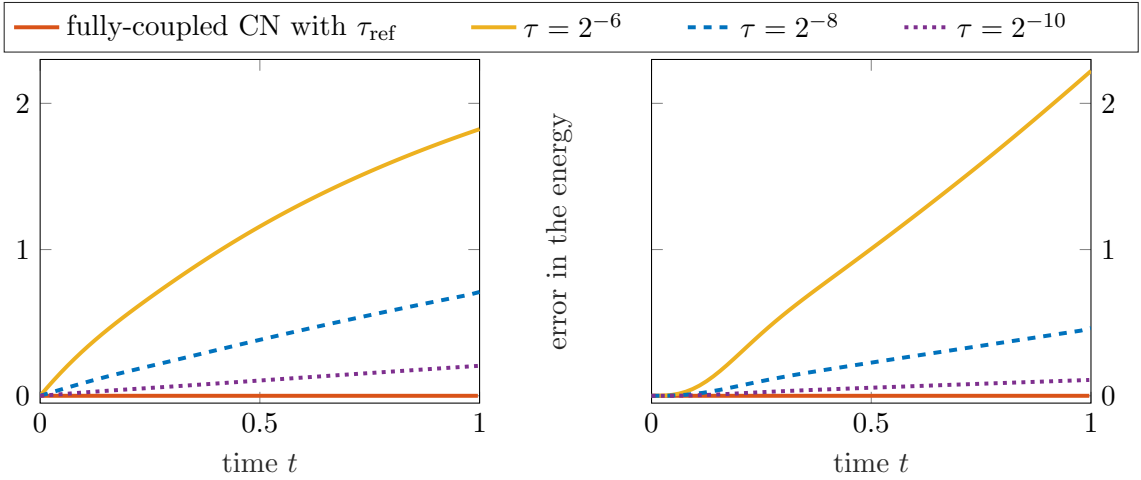
\begin{figure} 
%
%
\begin{tikzpicture}

\begin{axis}[%
width=2.3in,
height=1.8in,
scale only axis,
xmin=0,
xmax=1,
xtick={0, 0.5, 1},
xlabel style={font=\color{white!15!black}},
xlabel={time~$t$},
ymin=-0.05,
ymax=2.3,
yminorticks=true,
ylabel style={font=\color{white!15!black}},
legend columns = 4,
axis background/.style={fill=white},
legend style={at={(2.5,1.02)}, anchor=south east, legend cell align=left, align=left, draw=white!15!black}
]

\addplot [color=mycolor1, line width=1.5]
  table[row sep=crcr]{%
0	0\\
0.99609375	5.81756864903582e-14\\
};
\addlegendentry{fully-coupled CN with $\tau_\text{ref}$ \qquad}



\addplot [color=mycolor2, line width=1.5]
table[row sep=crcr]{%
0	0\\
0.015625	0.0548288328676998\\
0.03125	0.107557857611275\\
0.046875	0.158059745952927\\
0.0625	0.206428886756756\\
0.078125	0.252789643606862\\
0.09375	0.297266349809243\\
0.109375	0.339992252388512\\
0.125	0.381116624540268\\
0.140625	0.420803724152772\\
0.15625	0.459226242128518\\
0.171875	0.496556488597919\\
0.1875	0.532957609438724\\
0.203125	0.568576225203662\\
0.21875	0.603537244427131\\
0.234375	0.637941116129965\\
0.25	0.671863405043571\\
0.265625	0.705356304722105\\
0.28125	0.738451554023412\\
0.296875	0.771164180815935\\
0.3125	0.803496540448742\\
0.328125	0.835442216904661\\
0.34375	0.866989482592304\\
0.359375	0.898124143181091\\
0.375	0.928831708080136\\
0.390625	0.959098914102473\\
0.40625	0.988914685788026\\
0.421875	1.01827064282725\\
0.4375	1.04716126888656\\
0.453125	1.07558384449655\\
0.46875	1.10353822711526\\
0.484375	1.13102654032409\\
0.5	1.15805281565787\\
0.515625	1.18462261696973\\
0.53125	1.21074266874974\\
0.546875	1.2364205053811\\
0.5625	1.26166415617222\\
0.578125	1.2864818793307\\
0.59375	1.31088195546277\\
0.609375	1.3348725470042\\
0.625	1.3584616243015\\
0.640625	1.38165695257523\\
0.65625	1.40446612779168\\
0.671875	1.42689664467582\\
0.6875	1.44895597761175\\
0.703125	1.47065165544561\\
0.71875	1.49199131415684\\
0.734375	1.51298271642922\\
0.75	1.53363373342548\\
0.765625	1.55395229048959\\
0.78125	1.5739462840726\\
0.796875	1.59362348114245\\
0.8125	1.61299141427685\\
0.828125	1.63205728550616\\
0.84375	1.65082789005151\\
0.859375	1.66930956790366\\
0.875	1.6875081873302\\
0.890625	1.70542916049203\\
0.90625	1.72307748789947\\
0.921875	1.74045782578902\\
0.9375	1.75757456882121\\
0.953125	1.7744319397998\\
0.96875	1.79103407828137\\
0.984375	1.80738512080494\\
1	1.82348926681798\\
};
\addlegendentry{$\tau = 2^{-6}$\qquad}

\addplot [color=mycolor4, line width=1.5, dashed]
table[row sep=crcr]{%
0	0\\
0.0390625	0.0354120177805175\\
0.078125	0.0702137854926699\\
0.1171875	0.103479178613057\\
0.15625	0.134656662447472\\
0.1953125	0.163980904980292\\
0.234375	0.192214799663435\\
0.2734375	0.220126727119219\\
0.3125	0.248150645038267\\
0.3515625	0.276362916437189\\
0.390625	0.304638059429161\\
0.4296875	0.332806502077859\\
0.46875	0.360736226992717\\
0.5078125	0.388348922944587\\
0.546875	0.415607183505992\\
0.5859375	0.442497228232956\\
0.625	0.469015760032536\\
0.6640625	0.49516239319522\\
0.703125	0.520937238707685\\
0.7421875	0.546341849753255\\
0.78125	0.571380527325457\\
0.8203125	0.596059775908625\\
0.859375	0.620386161114606\\
0.8984375	0.644364460481678\\
0.9375	0.66799742269812\\
0.9765625	0.691286867709257\\
1	0.71\\
};
\addlegendentry{$\tau = 2^{-8}$\qquad}

\addplot [color=mycolor5, line width=1.5, dotted]
table[row sep=crcr]{%
0	0\\
0.009765625	0.00223296128211814\\
0.01953125	0.00447549533028546\\
0.029296875	0.00673099775475139\\
0.0390625	0.00899916335246242\\
0.048828125	0.0112682013404375\\
0.05859375	0.0135387139264656\\
0.068359375	0.0158044399891728\\
0.078125	0.0180566412021927\\
0.087890625	0.0202935267836728\\
0.09765625	0.0225090128636798\\
0.107421875	0.0246986164294127\\
0.1171875	0.0268597560406145\\
0.126953125	0.0289900253251871\\
0.13671875	0.0310883472963077\\
0.146484375	0.0331544892625666\\
0.15625	0.0351897977440276\\
0.166015625	0.0371960603353894\\
0.17578125	0.0391758556761155\\
0.185546875	0.0411327978785412\\
0.1953125	0.0430704705495875\\
0.205078125	0.044992831247106\\
0.21484375	0.0469040353801962\\
0.224609375	0.0488079194675888\\
0.234375	0.0507081491510566\\
0.244140625	0.0526080683974395\\
0.25390625	0.0545105206557395\\
0.263671875	0.0564178842322431\\
0.2734375	0.0583320538967844\\
0.283203125	0.0602544081553482\\
0.29296875	0.062185854449837\\
0.302734375	0.0641268920201168\\
0.3125	0.0660776445591629\\
0.322265625	0.0680379223578944\\
0.33203125	0.0700072943753711\\
0.341796875	0.0719851385994388\\
0.3515625	0.0739707026566792\\
0.361328125	0.0759631592529013\\
0.37109375	0.077961642837844\\
0.380859375	0.0799652845549366\\
0.390625	0.0819732416367871\\
0.400390625	0.0839847164642644\\
0.41015625	0.08599897106697\\
0.419921875	0.0880153360092946\\
0.4296875	0.0900332126338292\\
0.439453125	0.092052072794012\\
0.44921875	0.094071456939953\\
0.458984375	0.0960909695204637\\
0.46875	0.0981102731498815\\
0.478515625	0.100129082486769\\
0.48828125	0.102147157785373\\
0.498046875	0.104164298756344\\
0.5078125	0.106180338957594\\
0.517578125	0.108195140424848\\
0.52734375	0.110208588714982\\
0.537109375	0.112220588659594\\
0.546875	0.114231060829808\\
0.556640625	0.116239938566009\\
0.56640625	0.118247165392008\\
0.576171875	0.120252692738436\\
0.5859375	0.122256478045847\\
0.595703125	0.124258483271685\\
0.60546875	0.126258673733679\\
0.615234375	0.128257017229315\\
0.625	0.130253483357114\\
0.634765625	0.13224804297882\\
0.64453125	0.134240667870688\\
0.654296875	0.136231330616916\\
0.6640625	0.138220004633555\\
0.673828125	0.140206664172936\\
0.68359375	0.142191284356818\\
0.693359375	0.144173841381694\\
0.703125	0.146154312829462\\
0.712890625	0.148132677839211\\
0.72265625	0.150108917069826\\
0.732421875	0.152083012649961\\
0.7421875	0.154054948258639\\
0.751953125	0.156024709203832\\
0.76171875	0.157992282290361\\
0.771484375	0.159957655481691\\
0.78125	0.161920817553776\\
0.791015625	0.163881757873692\\
0.80078125	0.165840466238337\\
0.810546875	0.167796932642378\\
0.8203125	0.169751146965941\\
0.830078125	0.171703098699456\\
0.83984375	0.17365277679869\\
0.849609375	0.175600169647634\\
0.859375	0.177545265051713\\
0.869140625	0.17948805023378\\
0.87890625	0.181428511867631\\
0.888671875	0.183366636181337\\
0.8984375	0.185302409116282\\
0.908203125	0.187235816504372\\
0.91796875	0.18916684424105\\
0.927734375	0.191095478449348\\
0.9375	0.19302170563059\\
0.947265625	0.194945512793643\\
0.95703125	0.196866887562637\\
0.966796875	0.198785818274481\\
0.9765625	0.200702294077838\\
0.986328125	0.202616305031306\\
0.99609375	0.204527842184983\\
};
\addlegendentry{$\tau = 2^{-10}$\qquad} 

\end{axis}
%
%

\begin{axis}[%
width=2.3in,
height=1.8in,
at={(3.2in,0.in)},
scale only axis,
xmin=0,
xmax=1,
xtick={0, 0.5, 1},
xlabel style={font=\color{white!15!black}},
xlabel={time~$t$},
ymin=-0.05,
ymax=2.3,
yminorticks=true,
ylabel style={font=\color{white!15!black}},
yticklabel pos=right,
ylabel={error in the energy},
legend columns = 4,
legend style={at={(2.2,1.1)}, anchor=south east, legend cell align=left, align=left, draw=white!15!black}
axis background/.style={fill=white},
legend style={legend cell align=left, align=left, draw=white!15!black}
]

\addplot [color=mycolor1, line width=1.5]
  table[row sep=crcr]{%
0	0\\
0.99609375	5.81756864903582e-14\\
};



\addplot [color=mycolor2, line width=1.5]
table[row sep=crcr]{%
0	0\\
0.015625	0.000429698114170929\\
0.03125	0.00208646362099429\\
0.046875	0.00650579473763546\\
0.0625	0.0149167535955423\\
0.078125	0.0280379158030573\\
0.09375	0.0465627130513702\\
0.109375	0.0703411529489251\\
0.125	0.099421596491251\\
0.140625	0.132923295329529\\
0.15625	0.170394260233097\\
0.171875	0.210597596358614\\
0.1875	0.252904651660095\\
0.203125	0.296132528025601\\
0.21875	0.339768847928207\\
0.234375	0.382961151058149\\
0.25	0.425462863375834\\
0.265625	0.466823654282412\\
0.28125	0.5070642527345\\
0.296875	0.546051041997648\\
0.3125	0.583983510198876\\
0.328125	0.620899690245666\\
0.34375	0.657068150943628\\
0.359375	0.692576513236993\\
0.375	0.727676613953017\\
0.390625	0.762439237770022\\
0.40625	0.797057025181278\\
0.421875	0.831567531561365\\
0.4375	0.866096473305265\\
0.453125	0.900659120618402\\
0.46875	0.935329979137816\\
0.484375	0.97011439819874\\
0.5	1.00505567032442\\
0.515625	1.04015698198644\\
0.53125	1.07544232131887\\
0.546875	1.11091147747757\\
0.5625	1.14657453079856\\
0.578125	1.18242477991169\\
0.59375	1.21845861321148\\
0.609375	1.25466449024553\\
0.625	1.29103152121682\\
0.640625	1.32754790897404\\
0.65625	1.36420610763994\\
0.671875	1.40100418320683\\
0.6875	1.43794608574797\\
0.703125	1.47504540661853\\
0.71875	1.51231936758701\\
0.734375	1.5497927030539\\
0.75	1.58748815552525\\
0.765625	1.62542935794199\\
0.78125	1.66363214695984\\
0.796875	1.70210569124803\\
0.8125	1.74085002869054\\
0.828125	1.77985595393798\\
0.84375	1.81910652748051\\
0.859375	1.85858032126776\\
0.875	1.8982533179632\\
0.890625	1.93810314295117\\
0.90625	1.9781090805492\\
0.921875	2.01826052449921\\
0.9375	2.0585500976138\\
0.953125	2.09897983886693\\
0.96875	2.13956007139298\\
0.984375	2.18030779042128\\
1	2.22124465007392\\
};

\addplot [color=mycolor4, line width=1.5, dashed]
table[row sep=crcr]{%
0	0\\
0.0390625	0.000923742043299924\\
0.078125	0.00685521192587979\\
0.1171875	0.0204174102463832\\
0.15625	0.040932915361573\\
0.1953125	0.0653669261034051\\
0.234375	0.0904943980044832\\
0.2734375	0.114260176795417\\
0.3125	0.136063380208872\\
0.3515625	0.156237600724169\\
0.390625	0.175371781711876\\
0.4296875	0.19396175348104\\
0.46875	0.212285177653681\\
0.5078125	0.230484531795935\\
0.546875	0.248626045152677\\
0.5859375	0.266729462267668\\
0.625	0.284779093903796\\
0.6640625	0.302743763299488\\
0.703125	0.320620209013038\\
0.7421875	0.338450325911902\\
0.78125	0.356296974742119\\
0.8203125	0.374192124730831\\
0.859375	0.392112986992003\\
0.8984375	0.410004266116388\\
0.9375	0.427822211091017\\
0.9765625	0.445566012917938\\
1	0.463\\
};

\addplot [color=mycolor5, line width=1.5, dotted]
table[row sep=crcr]{%
0	0\\
0.009765625	4.14173325324896e-06\\
0.01953125	2.85038346317101e-05\\
0.029296875	9.61019998246826e-05\\
0.0390625	0.000226173863830947\\
0.048828125	0.000439237955874194\\
0.05859375	0.000745562411024192\\
0.068359375	0.00116225154118732\\
0.078125	0.00169713014815054\\
0.087890625	0.00235076928657163\\
0.09765625	0.0031336248602436\\
0.107421875	0.00403514752853917\\
0.1171875	0.00505386429377497\\
0.126953125	0.00618518807376445\\
0.13671875	0.0074103751804846\\
0.146484375	0.00872635580753567\\
0.15625	0.0101143344802259\\
0.166015625	0.0115609131094425\\
0.17578125	0.0130557418203137\\
0.185546875	0.0145793573330564\\
0.1953125	0.0161252585962948\\
0.205078125	0.0176768425201344\\
0.21484375	0.0192258585655849\\
0.224609375	0.0207656485614121\\
0.234375	0.0222834858031216\\
0.244140625	0.0237805613822792\\
0.25390625	0.0252482147119282\\
0.263671875	0.0266844644887203\\
0.2734375	0.0280910597361888\\
0.283203125	0.0294630499301909\\
0.29296875	0.0308050259604826\\
0.302734375	0.0321162125282437\\
0.3125	0.0333986623813796\\
0.322265625	0.0346556906403452\\
0.33203125	0.0358876555114169\\
0.341796875	0.037099502903772\\
0.3515625	0.0382914146303466\\
0.361328125	0.0394669346989382\\
0.37109375	0.0406286342690145\\
0.380859375	0.0417767534990281\\
0.390625	0.0429152800874935\\
0.400390625	0.0440441044465683\\
0.41015625	0.0451651290457282\\
0.419921875	0.0462799532749574\\
0.4296875	0.047388519518003\\
0.439453125	0.0484927139721933\\
0.44921875	0.0495923120921584\\
0.458984375	0.050688396890898\\
0.46875	0.0517813541228174\\
0.478515625	0.0528711377141073\\
0.48828125	0.05395882174802\\
0.498046875	0.0550438980953776\\
0.5078125	0.0561270532718212\\
0.517578125	0.0572084559043273\\
0.52734375	0.0582878866942429\\
0.537109375	0.0593659386370811\\
0.546875	0.0604422698285632\\
0.556640625	0.0615171085613473\\
0.56640625	0.0625904254816532\\
0.576171875	0.0636620514214785\\
0.5859375	0.0647321021244163\\
0.595703125	0.0658002376592601\\
0.60546875	0.0668665538579769\\
0.615234375	0.0679308009298993\\
0.625	0.0689928211540378\\
0.634765625	0.0700526979041443\\
0.64453125	0.0711101049503724\\
0.654296875	0.0721652006331879\\
0.6640625	0.0732179364157726\\
0.673828125	0.0742682440324631\\
0.68359375	0.0753164483026763\\
0.693359375	0.0763625338884437\\
0.703125	0.0774067205733546\\
0.712890625	0.0784492945747552\\
0.72265625	0.0794904115642407\\
0.732421875	0.0805304278845571\\
0.7421875	0.0815695270484116\\
0.751953125	0.0826079167198532\\
0.76171875	0.083645864438743\\
0.771484375	0.084683478539032\\
0.78125	0.0857208583619342\\
0.791015625	0.0867580675739119\\
0.80078125	0.0877951037098916\\
0.810546875	0.0888318748437764\\
0.8203125	0.0898682848272636\\
0.830078125	0.0909041853902512\\
0.83984375	0.0919393505004265\\
0.849609375	0.0929736341359066\\
0.859375	0.094006770550497\\
0.869140625	0.0950385343329625\\
0.87890625	0.0960688054666146\\
0.888671875	0.0970973210680688\\
0.8984375	0.0981239812247723\\
0.908203125	0.0991486961680961\\
0.91796875	0.100171344343285\\
0.927734375	0.101191998669653\\
0.9375	0.10221060441593\\
0.947265625	0.103227246460026\\
0.95703125	0.104242103929016\\
0.966796875	0.105255217014357\\
0.9765625	0.106266857328197\\
0.986328125	0.10727721306568\\
0.99609375	0.108286472601608\\
};

\end{axis}
\end{tikzpicture}%
	\caption{Development of the energy error over time for different step sizes~$\tau$ on a fixed spatial mesh. Plots show the results for the Lie--Euler splitting (left) and the Strang--CN splitting (right). }	
	\label{fig:kinetic:energy}	
\end{figure}

%
%
As mentioned in Remark~\ref{rem:kinetic:energy}, the original system has a Hamiltonian structure, which preserves the energy over time. We now investigate whether this is also true for the proposed splitting schemes. For this, we consider the fine triangulation, i.e., we consider the mesh parameter $h_\text{fine}$. 
Considering a Crank--Nicolson approximation of the fully-coupled system with mesh size $\tau_\text{ref}=2^{-12}$, the energy is preserved up to errors of the order $10^{-13}$. 
For the splitting schemes, on the other hand, the energy is not preserved, cf.~Figure~\ref{fig:kinetic:energy}. More precisely, the energy errors for both splitting approaches behave similar to the overall error, i.e., the errors are roughly given by $\mathcal{O}(\tau)$. 
%
\subsubsection*{A semi-linear example}
As a second experiment, we consider the semi-linear case. With the same parameters and initial data as in the linear setting, we now investigate the typical Allen--Cahn-type nonlinearity in the bulk, i.e., we define 
\[
	f_\Omega(t,u) = -u^3+u,\qquad
	f_\Gamma(t,p) = 0. 
\]
For both splitting schemes, we use an explicit treatment of the nonlinearity as explained in the previous subsections. This then leads to the same convergence results as in the linear case. Because of this, we omit the plot of the corresponding convergence history at this point and pass over to the more interesting case of a nonlinearity on the boundary. 
\subsubsection*{Comparison with other numerical schemes}
\begin{figure} 
%
%
\begin{tikzpicture}

\begin{axis}[%
width=3.8in, 
height=1.8in,
at={(2.45in,0.in)},
scale only axis,
xmode=log,
xmin=0.015,
xmax=810,
xminorticks=true,
xlabel style={font=\color{white!15!black}},
xlabel={computation time in~[$s$]},
ymode=log,
ymin=2e-04,
ymax=0.15,
yminorticks=true,
ylabel style={font=\color{white!15!black}},
ylabel={$L^\infty(L^2(\Omega))$-error in $u$},
axis background/.style={fill=white},
title style={font=\bfseries},
legend style={at={(1.45,0.35)}, anchor=south east, legend cell align=left, align=left, draw=white!15!black}
]


\addplot [color=mycolor1, line width=1.5]
table[row sep=crcr]{%
0.029320 4.6736555e-02\\
0.055596 1.3376470e-02\\
0.110877 3.6307275e-03\\
0.212154 9.8351388e-04\\
0.395246 2.4882771e-04\\
0.805043 6.1539768e-05\\
1.613355 1.4653456e-05\\
};
\addlegendentry{CN from \cite{HocL21}} 

\addplot [color=mycolor2, line width=1.5]
table[row sep=crcr]{%
0.030379 1.3088327e-01\\
0.060494 8.2938425e-02\\
0.113855 4.8394025e-02\\
0.229807 2.6558478e-02\\ 
0.457163 1.4005244e-02\\
0.847273 7.2129103e-03\\
1.663584 3.6654221e-03\\
3.339749 1.8487750e-03\\
6.548946 9.2863601e-04\\
12.994830 4.6541629e-04\\
26.242717 2.3298915e-04\\
%
};
\addlegendentry{IMEX Euler}

\addplot [color=mycolor4, line width=1.5]
table[row sep=crcr]{%
21.016373 1.2931427e-01\\
26.781638 8.2055574e-02\\
45.443415 4.7917061e-02\\
95.226497 2.6309871e-02\\
174.747421 1.3878680e-02\\
386.178592 7.1491810e-03\\
812.572360 3.6334825e-03\\
};
\addlegendentry{impl Euler}  
%
%
\addplot [color=mycolor2, line width=1.5, dashed]
table[row sep=crcr]{%
0.014109 6.8775994e-02\\
0.028566 4.4644727e-02\\
0.047547 2.7315095e-02\\
0.093733 1.6011500e-02\\
0.198417 9.0270620e-03\\
0.340964 4.9099868e-03\\
0.721221 2.5968858e-03\\
1.388249 1.3460115e-03\\
2.792629 6.8768683e-04\\
5.619392 3.4801901e-04\\ 	
12.629031 1.7512703e-04\\
%
};
\addlegendentry{Lie, IMEX Euler}

\addplot [color=mycolor4, line width=1.5, dashed]
table[row sep=crcr]{%
0.553210 6.7960813e-02\\
0.692906 4.4169142e-02\\
1.033493 2.7074224e-02\\
2.053059 1.5885544e-02\\
3.119480 8.9643572e-03\\
5.455487 4.8792854e-03\\
10.243355 2.5818812e-03\\
19.449359 1.3386487e-03\\
39.778645 6.8405226e-04\\
80.000711 3.4621537e-04\\	
154.562698 1.7422869e-04\\
%
};
\addlegendentry{Lie, impl Euler} 

\end{axis}
\end{tikzpicture}%
	\caption{Comparison of different solvers for kinetic boundary conditions with a nonlinearity on the boundary, namely the implicit--explicit Crank--Nicolson scheme, the implicit--explicit Euler scheme, the implicit Euler scheme, and two corresponding Lie splitting schemes. }
	\label{fig:kinetic:comparison}	
\end{figure}
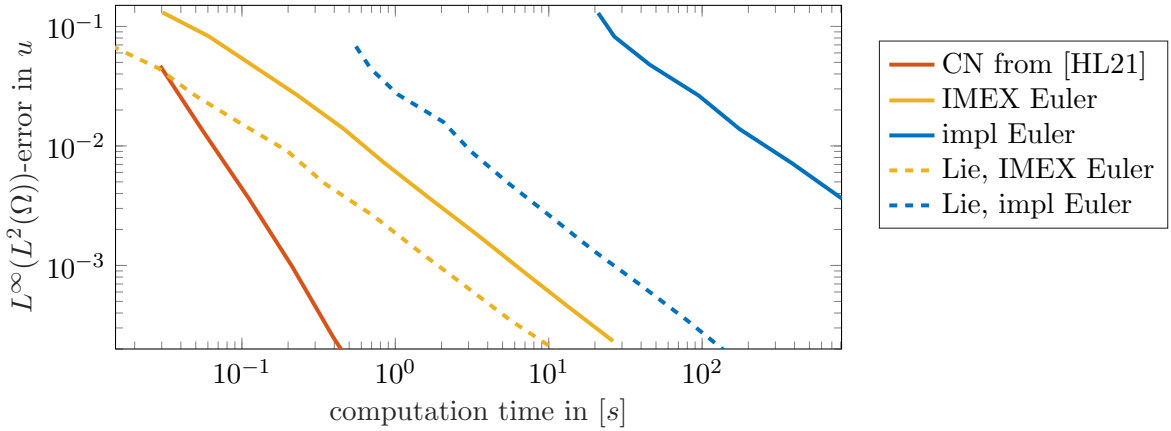

%
Finally, we would like to compare the proposed numerical schemes to standard methods applied to the fully-coupled system, i.e., to the original system without a splitting. Since bulk--surface splitting schemes are of particular interest when the solution on the boundary is strongly oscillatory or strongly nonlinear, we now investigate the case where the Allen--Cahn-type nonlinearity acts only on the boundary. Hence, we set  
\[
	f_\Omega(t,u) = 0,\qquad
	f_\Gamma(t,p) = -p^3+p. 
\]
All other parameters and the initial data remain unchanged. The following computations are performed on the coarse triangulation with mesh size $h_\text{coarse}$ and the errors are measured in the $L^\infty(L^2(\Omega))$-norm for $u$. The results in Figure~\ref{fig:kinetic:comparison} clearly show that the implicit-explicit Crank--Nicolson scheme from~\cite{HocL21} performs best, since it is the only scheme of second order. Moreover, no nonlinear system needs to be solved in this implicit--explicit version of the scheme. 
For the first-order schemes, we observe computational improvements due to the splitting. For the implicit--explicit Euler (where the nonlinearity is treated explicitly), the improvement of the splitting comes from the fact that we solve two smaller systems rather than a large one. For the implicit Euler scheme, the computational gain is even larger. Here, the splitting scheme exploits the fact that the nonlinearity only acts on the boundary. 
%

As a summary, the bulk--surface splitting schemes for the wave equation with kinetic boundary conditions proposed in this paper only reach first-order convergence. Whereas the Lie--Euler scheme shows expectable convergence results (independent of the spatial mesh size), the convergence order of the Strang--CN scheme is reduced and $h$-dependent. In the subsequent section, we will see that acoustic boundary conditions allow for second-order splitting schemes. 
%
%
\section{Splitting Schemes for Acoustic Boundary Conditions}
\label{sect:schemesAcoustic}
This section is devoted to the design of splitting methods for the wave equation with acoustic boundary conditions~\eqref{eq:PDAE:acoustic}. In contrast to the previous section, a spatial discretization by finite elements yields an ordinary differential equation rather than a DAE. 
We introduce a splitting of bulk and surface dynamics and combine this with suitable time stepping schemes. Here, we observe convergence rates of order one for Lie and order two for Strang splitting. 
%
%
\subsection{Semi-discrete system and bulk--surface splitting}
For the spatial discretization, we consider again bulk--surface finite elements as introduced in Section~\ref{sect:schemesKinetic:FEM}. The semi-discrete version of~\eqref{eq:PDAE:acoustic} then reads 
\begin{align}
\label{eq:semidisc:acoustic}
	\begin{bmatrix} M_\Omega &  \\  & M_\Gamma \end{bmatrix}
	\begin{bmatrix} \ddot u \\ \ddot \delta  \end{bmatrix}
	+ \begin{bmatrix} & -B^T \\ B & \end{bmatrix}
	\begin{bmatrix} \dot u \\ \dot \delta  \end{bmatrix}
	+ \begin{bmatrix} A_\Omega &  \\  & A_\Gamma \end{bmatrix}
	\begin{bmatrix} u \\ \delta \end{bmatrix} 
	&= \begin{bmatrix} f_\Omega(u) \\ f_\Gamma(\delta) \end{bmatrix}.
\end{align}
Here, we use the same notion as before with the exception of $B$, which (assuming once more the specific ordering of the basis functions) has the form~$B=[\, 0\ \ M_\Gamma]$. Note that~\eqref{eq:semidisc:acoustic} is an ordinary differential equation such that no regularization as for kinetic boundary conditions is necessary. 

For a splitting of~\eqref{eq:semidisc:acoustic} into a bulk and a boundary problem, we consider the first time interval~$[0,\tau]$ with initial data~$u(0) = u^0$, $\dot u(0) = \dot u^0$, $\delta(0) = \delta^0$, and~$\dot \delta(0) = \dot \delta^0$.
We would like to emphasize that the variable $u$ is not decomposed, i.e., the bulk problem also contains the boundary. The following splitting approach aims to decouple the dynamics of $u$ and $\delta$ on the boundary. 

{\em Subsystem 1 (bulk problem)}: 
In the bulk problem, we solve  
\begin{align}
\label{eq:acoustic:bulk}
  M_\Omega{\ddot u} + A_\Omega u 
  = f_\Omega(u) + B^T \dot \delta^0 
\end{align}
with initial conditions $u(0) = u^0$ and~$\dot u(0) = \dot u^0$. Within this subproblem, the displacement of the boundary in normal direction (i.e., $\delta$) remains unchanged. The solution of~\eqref{eq:acoustic:bulk} yields the approximation~$u(\tau)$ as well as~$\dot u(\tau)$. 

{\em Subsystem 2 (boundary problem)}: 
In the boundary problem, $u$ remains unchanged and we solve
\begin{align}
\label{eq:acoustic:boundary}
  M_\Gamma \ddot \delta + A_\Gamma \delta
  = f_\Gamma(\delta) - B \dot u(\tau)
\end{align}
with initial conditions $\delta(0) = \delta^0$ and~$\dot \delta(0) = \dot \delta^0$, leading to the new approximations~$\delta(\tau)$ and $\dot \delta(\tau)$. 
%
%
\subsection{Fully-discrete splitting schemes}\label{sect:schemesAcoustic:fullyDiscrete}
To obtain practical applicable bulk--surface splitting methods, we further need to discretize the subsystems in time. As for kinetic boundary conditions, we pursue two possibilities: a Lie splitting approach together with an implicit Euler discretization and a Strang splitting with a Crank--Nicolson discretization. 
%
%
\subsubsection{Lie splitting with implicit Euler discretization}\label{sect:schemesAcoustic:fullyDiscrete:Euler}
We consider a Lie splitting solving on each subinterval of length $\tau$ first the bulk problem and afterwards the boundary problem. Both subsystems are discretized by the implicit Euler scheme~\eqref{eq:implicitEuler}, where the possible nonlinearity is treated explicitly. Hence, given the approximations~$u^n$, $w^n\coloneqq\dot u^n$, $\delta^n$, and~$\zeta^n\coloneqq\dot \delta^n$ at time $t^n$, the first subsystem reads 
\begin{align*}
  (M_\Omega + \tau^2 A_\Omega)\, w^{n+1}  
  &= M_\Omega w^n - \tau A_\Omega u^n + \tau\, (f_\Omega^{n+1} + B^T\zeta^n), \\
  u^{n+1}
  &= u^n + \tau\, w^{n+1}
\end{align*}
with $f_\Omega^{n+1} = f_\Omega(t^{n+1},u^n)$ in the semi-linear case. This then gives new approximations $u^{n+1}$ and $w^{n+1}$. 
The discretization of the second subsystem yields 
\begin{align*}
  (M_\Gamma + \tau^2 A_\Gamma)\, \zeta^{n+1}  
  &= M_\Gamma \zeta^n - \tau A_\Gamma \delta^n + \tau\, (f_\Gamma^{n+1} - B w^{n+1}), \\
  \delta^{n+1}
  &= \delta^{n} + \tau\, \zeta^{n+1}.
\end{align*}
and provides updates $\delta^{n+1}$ and $\zeta^{n+1}$ at time $t^{n+1}$. 
%
%
\subsubsection{Strang splitting with Crank--Nicolson discretization}\label{sect:schemesAcoustic:fullyDiscrete:CN}
For the symmetric splitting, we first solve the bulk problem~\eqref{eq:acoustic:bulk} on the half interval, then the boundary problem~\eqref{eq:acoustic:boundary} on the entire interval, and finally the remaining part of the bulk problem. Here, all subsystems are discretized by the implicit--explicit Crank--Nicolson scheme presented in~\eqref{eq:CrankNicolson}. 

Given the approximations~$u^n$, $w^n\coloneqq\dot u^n$, $\delta^n$, and~$\zeta^n\coloneqq\dot \delta^n$ at time $t^n$, we first compute 
\begin{align*}
  (M_\Omega + \tfrac{\tau^2}{16} A_\Omega)\, w^{n+1/4}  
  &= M_\Omega w^n - \tfrac{\tau}{4} A_\Omega u^n + \tfrac{\tau}{4}\, (f_\Omega^{n} + B^T\zeta^n), \\
  u^{n+1/2}
  &= u^n + \tfrac{\tau}{2}\, w^{n+1/4}, \\
  M_\Omega w^{n+1/2}  
  &= 2M_\Omega w^{n+1/4} - M_\Omega w^n + \tfrac{\tau}{4}\, (f_\Omega^{n+1/2} - f_\Omega^{n}).
\end{align*}
Recall that, in the semi-linear case, $f_\Omega^{n+1/2}$ includes the updated state~$u^{n+1/2}$. With the obtained approximation~$w^{n+1/2}$, we may then compute, as discretization of the second subsystem, 
\begin{align*}
  ( M_\Gamma + \tfrac{\tau^2}{4} A_\Gamma )\, \zeta^{n+1/2} 
  &= M_\Gamma \zeta^n - \tfrac{\tau}{2} A_\Gamma \delta^n + \tfrac{\tau}{2} ( f_\Gamma^n - B w^{n+1/2} ), \\
  \delta^{n+1}
  &= \delta^n + \tau\, \zeta^{n+1/2}, \\ 
  M_\Gamma \zeta^{n+1} 
  &= 2M_\Gamma \zeta^{n+1/2} - M_\Gamma \zeta^n + \tfrac{\tau}{2}\, ( f_\Gamma^{n+1}-f_\Gamma^n ),
\end{align*}
which yields updates $\delta^{n+1}$ and $\zeta^{n+1}$ at time $t^{n+1}$. 
Finally, the third subsystem reads 
\begin{align*}
  (M_\Omega + \tfrac{\tau^2}{16} A_\Omega)\, w^{n+3/4}  
  &= M_\Omega w^{n+1/2} - \tfrac{\tau}{4} A_\Omega u^{n+1/2} + \tfrac{\tau}{4}\, (f_\Omega^{n+1/2} + B^T\zeta^{n+1}), \\
  u^{n+1}
  &= u^{n+1/2} + \tfrac{\tau}{2}\, w^{n+3/4}, \\
  M_\Omega w^{n+1}  
  &= 2M_\Omega w^{n+3/4} - M_\Omega w^{n+1/2} + \tfrac{\tau}{4}\, (f_\Omega^{n+1} - f_\Omega^{n+1/2})
\end{align*}
and provides updates $u^{n+1}$ and $w^{n+1}$. 

In the following experiments, we show that this approach indeed yields a second-order method. 
%
%
\subsection{Numerical experiment}\label{sect:schemesAcoustic:numerics}
In this final section on acoustic boundary conditions, we provide a numerical experiment, which indicates convergence rates of order one and two, respectively. 
As in Section~\ref{sect:schemesKinetic:numerics}, the computational domain is given by the unit disc and we compare the numerical solutions to a reference solution resulting from a Crank--Nicolson discretization of the fully-coupled system~\eqref{eq:semidisc:acoustic} with step size $\tau_\text{ref}=2^{-12}$. 
\subsubsection*{Initial data and parameters}
Following \cite[Sect 8.4]{HipK20}, we set as initial data
\[
	u^0(x,y) = 0, \qquad
	w^0(x,y) = 2\pi\, (x^2+y^2)^{0.6} 
\]
for the acoustic velocity potential and 
\[
	\delta^0(x,y) = \tfrac{k}{2\pi}\, (x^2+y^2)^{0.6}, \qquad
	\zeta^0(x,y) = 0
\]
for the displacement of the boundary in normal direction. For the remaining parameters in~\eqref{eq:acousticBC}, we set again 
\[
	\beta = 1, \qquad
	\kappa = 1, \qquad 
	T = 1.
\]
Mesh parameters and step sizes are chosen as in the previous section, i.e., we consider two triangulations with $h_\text{coarse} \approx 0.094$ and $h_\text{fine} \approx 0.025$, respectively, and the step sizes vary from $2^{-4}$ to $2^{-11}$. 
\subsubsection*{Order of convergence for a semi-linear example}
\begin{figure} 
%
%
\begin{tikzpicture}
\begin{axis}[%
width=2.3in,
height=1.8in,
scale only axis,
xmode=log,
xmin=0.00045,
xmax=0.0625,
xminorticks=true,
xlabel style={font=\color{white!15!black}},
xlabel={step size~$\tau$},
ymode=log,
ymin=2e-05,
ymax=4e-00,
yminorticks=true,
ylabel style={font=\color{white!15!black}},
axis background/.style={fill=white},
title style={font=\bfseries},
legend columns = 2,
legend style={at={(1.89,1.02)}, anchor=south east, legend cell align=left, align=left, draw=white!15!black}
]


\addplot [color=mycolor2, line width=1.5]
table[row sep=crcr]{%
	0.0625	2.3850989\\
	0.03125	1.5951369\\
	0.015625	0.99282821\\
	0.0078125	0.58467114\\
	0.00390625	0.32913945\\
	0.001953125	0.17867937\\
	0.0009765625	0.094344504\\
	0.00048828125	0.048842457\\
};
\addlegendentry{Lie--Euler with $h_\text{fine}$\qquad} 

\addplot [color=mycolor2, line width=1.5, dashed]
table[row sep=crcr]{%
	0.0625	2.3215993\\
	0.03125	1.5380265\\
	0.015625	0.94106449\\
	0.0078125	0.54199607\\
	0.00390625	0.29909716\\
	0.001953125	0.16033279\\
	0.0009765625	0.084149408\\
	0.00048828125	0.043425268\\
};
\addlegendentry{Lie--Euler with $h_\text{coarse}$} 


\addplot [color=mycolor5, line width=1.5]
table[row sep=crcr]{%
	0.0625	1.1193772\\
	0.03125	0.26085509\\
	0.015625	0.065502609\\
	0.0078125	0.016770439\\
	0.00390625	0.0043791645\\
	0.001953125	0.0011090455\\
	0.0009765625	0.00026842358\\
	0.00048828125	6.2725744e-05\\
};
\addlegendentry{Strang--CN with $h_\text{fine}$\qquad} 

\addplot [color=mycolor5, line width=1.5, dashed]
table[row sep=crcr]{%
	0.0625	0.87559264\\
	0.03125	0.22321237\\
	0.015625	0.056393137\\
	0.0078125	0.0141267\\
	0.00390625	0.0035286513\\
	0.001953125	0.00087796684\\
	0.0009765625	0.00021566893\\
	0.00048828125	5.1988583e-05\\
};
\addlegendentry{Strang--CN with $h_\text{coarse}$} 


\addplot [color=gray, dotted, line width=1.0pt]
  table[row sep=crcr]{%
	0.0625	0.9\\
	0.0000625	0.0009\\
};

\addplot [color=gray, dotted, line width=1.0pt]
  table[row sep=crcr]{%
  	0.0625	0.1\\
  	0.000625	0.00001\\
};

\end{axis}

%
%
\begin{axis}[%
width=2.3in,
height=1.8in,
at={(2.45in,0.in)},
scale only axis,
xmode=log,
xmin=0.00045,
xmax=0.06,
xminorticks=true,
xlabel style={font=\color{white!15!black}},
xlabel={step size~$\tau$},
ymode=log,
ymin=2e-05,
ymax=4e-00,
yminorticks=true,
yticklabel pos=right,
ylabel style={font=\color{white!15!black}},
axis background/.style={fill=white},
title style={font=\bfseries},
]


\addplot [color=mycolor2, line width=1.5]
table[row sep=crcr]{%
	0.0625	2.1260983\\
	0.03125	0.82429221\\
	0.015625	0.36045816\\
	0.0078125	0.16781944\\
	0.00390625	0.08085417\\
	0.001953125	0.039668611\\
	0.0009765625	0.019648588\\
	0.00048828125	0.0097803996\\
};

\addplot [color=mycolor2, line width=1.5, dashed]
table[row sep=crcr]{%
	0.0625	2.1080744\\
	0.03125	0.81545273\\
	0.015625	0.35632267\\
	0.0078125	0.16587696\\
	0.00390625	0.079922775\\
	0.001953125	0.039214168\\
	0.0009765625	0.019421106\\
	0.00048828125	0.0096641696\\
};


\addplot [color=mycolor5, line width=1.5]
table[row sep=crcr]{%
	0.0625	1.1696153\\
	0.03125	0.24185029\\
	0.015625	0.059901931\\
	0.0078125	0.014960433\\
	0.00390625	0.0037396333\\
	0.001953125	0.00093473878\\
	0.0009765625	0.00023380876\\
	0.00048828125	5.8545814e-05\\
};

\addplot [color=mycolor5, line width=1.5, dashed]
table[row sep=crcr]{%
	0.0625	0.31480278\\
	0.03125	0.08004494\\
	0.015625	0.02019255\\
	0.0078125	0.0050564286\\
	0.00390625	0.0012644082\\
	0.001953125	0.00031615723\\
	0.0009765625	7.9084237e-05\\
	0.00048828125	1.9818589e-05\\
};


\addplot [color=gray, dotted, line width=1.0pt]
table[row sep=crcr]{%
	0.0625	0.9\\
	0.0000625	0.0009\\
};

\addplot [color=gray, dotted, line width=1.0pt]
table[row sep=crcr]{%
	0.0625	0.1\\
	0.000625	0.00001\\
};

\end{axis}
\end{tikzpicture}%
\caption{Convergence history for the wave equation with acoustic boundary conditions and Allen--Cahn-type nonlinearity on the boundary. Plots show the $L^\infty(L^2(\Omega))$ error of $u$ (left) and the $L^\infty(L^2(\Gamma))$ error of~$\delta$ (right) for the two splitting schemes and two spatial mesh sizes $h_\text{fine}$ (solid) and $h_\text{coarse}$ (dashed). The dotted lines indicate orders 1 and 2, respectively.}
\label{fig:acoustic:nonlinear}
\end{figure}
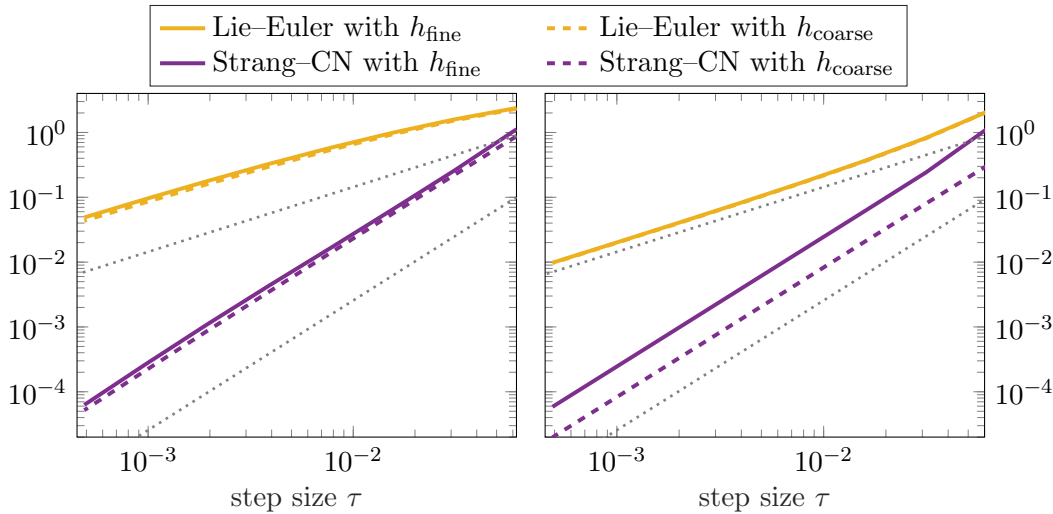

For the numerical investigation of the convergence orders of the proposed bulk--surface splitting schemes, we consider a semi-linear example with a nonlinearity of Allen--Cahn-type on the boundary. For this, we define the right-hand sides as 
\[
	f_\Omega(t,u) = 0,\qquad
	f_\Gamma(t,\delta) = -\delta^3+\delta. 
\]
The resulting convergence history is presented in Figure~\ref{fig:acoustic:nonlinear} and shows first-order convergence for the Lie--Euler scheme introduced in Section~\ref{sect:schemesAcoustic:fullyDiscrete:Euler} and second-order convergence for the Strang--CN scheme from Section~\ref{sect:schemesAcoustic:fullyDiscrete:CN}. 
%
%
%

Apart from the order of convergence, we are once more interested in a possible $h$-dependence of the proposed methods. For the Lie--Euler scheme, the dashed and solid lines in Figure~\ref{fig:acoustic:nonlinear} are on top of each other, indicating that the convergence is independent of the chosen spatial mesh size. Hence, the method is stable in this regard and it seems reasonable to prove first-order convergence of the Lie splitting, maybe even on operator level. This, however, is subject of future research. The Strang--CN splitting, on the other hand, shows a slight~$h$-dependence, at least for the displacement of the boundary in normal direction~$\delta$. Recall that this means that the error grows if the spatial mesh is refined. As a consequence, a spatial refinement calls for a smaller step size as well. Moreover, the Strang scheme can only be stable if some kind of a CFL-type condition holds and a second-order convergence result cannot count independently of the mesh size (and hence, not on operator level). 
\begin{remark} 
As one may expect, the convergence order drops to one if we discretize the subsystems of the Strang splitting with an implicit or implicit--explicit Euler scheme. 
\end{remark}
\subsubsection*{Conservation of energy}
%
%
Finally, we would like to comment on the evolution of the energy introduced in Remark~\ref{rem:acoustic:energy}. Similar as for kinetic boundary conditions, the energy is not preserved for the proposed splitting schemes but behaves as the overall error. This means that the energy error is of order $\mathcal{O}(\tau)$ for the Lie--Euler scheme and $\mathcal{O}(\tau^2)$ for the Strang--NC scheme. For the latter, however, the error grows if the spatial mesh is refined. More precisely, going from $h_\text{coarse}$ to $h_\text{fine}$, which correlate by a factor of around $4$, the error in the energy grows by a factor of $4$ as well.  
%
%
\section{Conclusion}\label{sect:conclusion}
Within this paper, we have derived fully-discrete splitting schemes for the wave equation with (non-local) kinetic and acoustic boundary conditions. 
In the case of acoustic boundary conditions, Lie and Strang splitting show the expected convergence rates if the subsystems are discretized accordingly. 
For kinetic boundary conditions, we considered a reformulation as PDAE for the construction of splitting schemes. In this case, the proposed schemes do not exceed order one. 

Future research calls for a rigorous error analysis of the proposed methods to confirm the observed convergence rates. Moreover, the construction of second-order splitting schemes for kinetic boundary conditions remains an open problem. But also for acoustic boundary conditions further improvements are desirable, namely a second-order scheme which does not depend on the spatial mesh size.
%
%
\bibliographystyle{alpha}
\bibliography{bib_dynBC}
\end{document}